\documentclass{article}  
\usepackage[top=15mm, bottom=25mm]{geometry}
\usepackage{authblk}
\usepackage{graphicx}
\usepackage[usenames,dvipsnames]{xcolor}
\usepackage{tikz}
\usepackage{booktabs}
\usepackage{array}
\usepackage{paralist}
\usepackage{natbib}
\usepackage{hyperref}
\usepackage{amsmath, amssymb}
\newcommand{\st}{\text{subject to}}
\newcommand{\CI}{\operatorname{CI}}

\begin{document}
  
\title{Rank-optimal weighting or ``How to be best in the OECD Better Life Index?''  
}

\author[*]{Jan Lorenz}
\author[**]{Christoph Brauer}
\author[**]{Dirk A. Lorenz}
\affil[*]{Jacobs University Bremen, Focus Area Diversity,  
              Campus Ring 1, 28759 Bremen, Germany \\
              \texttt{post@janlo.de}} 
\affil[**]{Technische Universit\"at Braunschweig, Institut f\"ur Analysis und Algebra, 38106 Braunschweig, Germany}

\date{Aug 15, 2016}

\maketitle

\begin{abstract}
We present a method of rank-optimal weighting which can be used to explore the best possible position of a subject in a ranking based on a composite indicator by means of a mathematical optimization problem. As an example, we explore the dataset of the OECD Better Life Index and compute for each country a weight vector which brings it as far up in the ranking as possible with the greatest advance of the immediate rivals. The method is able to answer the question ``What is the best possible rank a country can achieve with a given set of weighted indicators?'' Typically, weights in composite indicators are justified normatively and not empirically. Our approach helps to give bounds on what is achievable by such normative judgments from a purely output-oriented and strongly competitive perspective. The method can serve as a basis for exact bounds in sensitivity analysis focused on ranking positions. 

In the OECD Better Life Index data we find that 19 out the 36 countries in the OECD Better Life Index 2014 can be brought to the top of the ranking by specific weights. We give a table of weights for each country which brings it to its highest possible position. Many countries achieve their best rank by focusing on their strong dimensions and setting the weights of many others to zero. Although setting dimensions to zero is possible in the OECD's online tool, this contradicts the idea of better life being multidimensional in essence. We discuss modifications of the optimization problem which could take this into account, e.g. by allowing only a minimal weight of one. 

Methods to find rank-optimal weights can be useful for various multidimensional datasets like the ones used to rank universities or employers.

\end{abstract}


\section{Introduction}

Composite social indicators are increasingly used to assess and rank countries in public and scientific debates. Almost all composite indicators are weighted averages of a number of dimensions, where weights are set ad hoc with a brief normative justification. Typically, weights are set equal. The ranks of countries is what is mostly looked at, despite all warnings that weights are somehow arbitrarily set. Readers of rankings might know that small variations in weights can cause large shifts in ranks, but usually it is unknown how large the impact of the choice of weights can be. Consequently, this question is of particular interest for sensitivity analysis as well as for the public debate. The problem of rank-optimal weighting is to find weights which bring a particular country as far up in the ranking as possible. The problem is surprisingly complicated, but feasible. We present a method to do this and take the OECD Better Life Index dataset as an example. 

The OECD Better Life Index \citep{OECD2013BetterLifeIndex} is a collection of eleven dimensions measuring essential indicators to well-being, from health and education to local environment, personal security, and overall satisfaction with life, as well as more traditional measures such as income. Based on the recommendation of the Commission on the Measurement of Economic Performance and Social Progress \citep{Stiglitz.Sen.ea2009measurementeconomicperformance} the OECD Better Life Initiative sees well-being as a multidimensional concept. Consequently, they offer a tool on their website where everyone can build his or her own Better Life Index by assigning weights to the dimensions of well-being. This interactive tool is to engage people in the debate on well-being and, through this process, to learn what matters the most to them. The tool invites to compare well-being across countries according to the weights given to the dimensions: Housing, Income, Jobs, Community, Education, Environment, Civic Engagement, Health, Life Satisfaction, Safety, and Work-Life Balance.
Users are also invited to share their weights with others. 

This interactive approach has been adopted by the Legatum Prosperity Index \citep{LIF2014LegatumProsperityIndex} which measures prosperity of countries by eight sub-indices. 
In their report they state ``We offer you, the reader, the opportunity to assign your own weightings to each of the sub-indices, and to see how the rankings change according to these different weightings.'' 
Although it is acknowledged that the prosperity of countries in different states of development might be measured differently it is seen more crucial to measure each country by the same yard stick, because different weights for each country would make country rankings incomparable. This is an argument for an equal weight in each dimension for all countries, but not for equal weights across all dimensions. 

Other multidimensional societal concepts measured similarly are social cohesion, social progress, or human development. Social cohesion is measured in nine dimensions in a study called ``Cohesion Radar'' by the Bertelsmann Foundation \citep{Dragolov.Ignacz.ea2013InternationalComparisonof}. Social progress is measured in twelve dimensions by the Social Progress Index \citep{SPI2014SocialProgressIndex}. Both measurements concepts use simple averaging of dimensions with equal weights to compute their final scores. Analogously, the Human Development Index \citep{UNDP2014HumanDevelopmentIndex} is the simple geometric mean of three dimensions with equal weights. 

Common to the measurement concepts of better life, prosperity, social cohesion, social progress, human development and many other composite indices is that their final aggregation of dimension to a one-dimensional scale follows a formative composition paradigm and not a reflective one \citep{Bollen.Lennox1991Conventionalwisdommeasurement}. That means all dimensions together form the construct instead of the construct being reflected in every single dimension. Thus, the construction of the measurement is theoretically postulated and not empirically validated. Thus the door is open for normative discussions about the relative weights of dimensions. Different weights lead to different rankings and none has validated priority unless users achieve consensus on its theoretical grounding. When composite indicators are used to assess the quality of policies, governments or political systems such consensus is unlikely if not impossible. The open approach of the OECD Better Life Index is a way to avoid fixing weights and consequently avoiding the focus on one ranking only. 

Despite of the arbitrariness of normative postulations, composite indicators are a well received tool, as \citet{Saisana.Saltelli.ea2005UncertaintyandSensitivity} state: ``[...] the temptation of stakeholders and practitioners to summarize complex and sometime elusive processes (e.g. sustainability, single market policy, etc.) into a single figure to benchmark country performance for policy consumption seems likewise irresistible.''  Also others speak in favor of composite indicators \citep{Booysen2002OverviewandEvaluation,Saltelli2007CompositeIndicatorsbetween}. The OECD Handbook for constructing composite indicators \citep{Nardo.Saisana.ea2005HandbookConstructingComposite} states that composite indicators can summarize multi-dimensional issues in view of supporting decision-makers. They can also make interpretation easier compared to the interpretation of various indicators in parallel. Such indicators are and will be used for ranking countries for benchmarking exercises and for the assessment of progress over time. In this way, composite indicators trigger societal debates on complex issues in the general public e.g. through simplified narratives in the media.

The OECD Handbook further states that the relative importance of dimensions will always remain a source of contention because the assignment of weights is essentially a value judgment no matter what method is used. 

In summary, composite indicators will be used especially in public and political communication, but the weighting of the dimensions remains a value judgment. In the political realm, this raises the question to what extent such a value judgment can serve the simple purpose of bringing a particular country to the top of the list or at least as high up as possible. Readers of analyses based on composite indicator rankings might speculate about the leverage that different weights have on the ranking. Politicians of poorly performing countries might even blame the choice of weights to be the main reason of poor performance. Our method easily computes the precise upper bounds for the best rank some country can achieve.

In the following, we formulate the rank optimal weighting problem as a precise mathematical optimization problem and use it to compute rank-optimal weights for all countries in the OECD Better Life dataset. With these results and some examples we discuss potential uses of the method.

\section{Method}

Following the OECD Handbook \citep{Nardo.Saisana.ea2005HandbookConstructingComposite}, let us consider a dataset of $Q$ indicators and $C$ countries and denote 
$I_{q,c}$ to be the normalized value of sub-indicator (or dimension) $q \in \{1,\dots,Q\}$ for country $c \in \{1,\dots,C\}$. 
For weights $w_1,\dots,w_Q \geq 0$ the \emph{Composite Index} for country $c$ is defined as
\[
\CI_c = \sum_{q=1}^Q w_q I_{q,c}.
\]

How can we formulate the problem of rank-optimal weighting in precise mathematical terms? 
Let us first discuss some simple special cases: If country $c$ has the highest value in dimension $q$ ($I_{q,c} = \max_{c}I_{q,c}$) then a solution with $w_q=1$ and all other weights equal zero brings country $c$ naturally to the top of the list. Thus, the best countries in one of the dimensions naturally have an easy solution of weights which makes them lead the ranking. But also countries with no top score in any dimension might be able to reach the top of the list. As an example consider a country which has second best rankings in two dimensions. When the top countries in these two dimension are different and each of them has low scores on the other dimensions, then equal weights for these two dimensions (and zero to the others) would bring the country above the other two in the list. The problem of rank optimal weighting for a particular country is thus not always simple to solve by hand.

In the following we are going to describe a computational model to find the highest rank a country can achieve (and also to find the corresponding weights). To that end, we are going to formulate the problem as a problem of mathematical optimization, i.e. a problem of maximizing a certain objective function subject to certain constraints for the corresponding optimization variables.

\subsection{First order rank maximization for country $c$} \label{subsec:rankmaximization}

Let us fix some country $c\in\{1,\dots,C\}$. We aim at finding the weights such that the country $c$ has the highest possible rank. In other words: country $c$ dominates as much countries as possible in the ranking. We denote this number, i.e. the number of countries that have a lower (or equal) score than country $c$, by $R_c$ and try to maximize it over all possible weights. Consequently, $R_c$ is the number of country indices $k\in\{1,\dots,C\}\setminus\{c\}$ for which it holds that $\CI_{c}\geq \CI_{k}$. For the feasible set of weights we assume for now that we can choose any non-negative weight. We add the constraint that the weights should sum up to ten to avoid that all countries have the same score $\CI=0$ for $w_q=0$, i.e. we demand $\sum_q w_q = 10$ and $w_q\geq 0$. Note that we could also demand that the $w_q$ sum to any number different from ten since this would only scale the score of all countries in the same way.\footnote{The choice $\sum_q w_q = 10$ corresponds to the scaling in OECD Better Life Index since therein each sub-indicator is normalized such that $0\leq I_{q, c}\leq 1$ and the composite index for each country can take values $0\leq \CI_{c} \leq 10$.} To write down an optimization problem we need to express our objective $R_c$ in terms of the vector of optimization variables $\vec{w} = [w_1,\dots,w_Q]$. In principle we could write
\[
R_c = \#\{k\ :\ \CI_c\geq \CI_k\}
\]
where $\#$ stands for ``number of elements in the set''. Note that this expresses $R_c$ in terms of the weights $w_q$ since $\CI_c\geq \CI_k$ if and only of $\sum_{q = 1}^{Q}I_{q, c}w_{q}\geq \sum_{q = 1}^{Q}I_{q, k}w_{q}$. This gives rise to the optimization problem
\begin{align*}
    R_c^{*} & = \max\limits_{\vec{w}} \#\Big\{k\ :\ \sum_{q = 1}^{Q}I_{q, c}w_{q}\geq \sum_{q = 1}^{Q}I_{q, k}w_{q}\Big\}\\
  & \begin{aligned}
    \text{subject to} &&\sum\limits_{q=1}^Q w_q&= 10\\
    &&w_q&\geq 0, &q=1,\dots, Q.
  \end{aligned}
\end{align*}
In this form, the objective function is not a nice function of the optimization variable $\vec{w}$ in the sense that it is not a continuous or differentiable function. In principle one could still try to solve the problem by ``zeroth order methods'' which try to explore the search space for $\vec{w}$ in a systematic fashion (e.g. Nelder-Mead simplex search, simulated annealing or the like). However, these methods usually become very slow even when the number of variables (i.e. the number $Q$) is mildly large. Much more powerful algorithms exist in the case where the objective function is differentiable. To convert the problem into one which is more amenable for optimization algorithms we employ two tricks: \emph{Auxiliary variables} and \emph{relaxation}. First we introduce additional variables $\vec{z} = [z_k]_{k\in\{1,\dots,C\}\setminus\{c\}}$. The variable $z_k$ shall indicate if the score of country $c$ is at least equal to the score of country $k$, i.e. we want
\[
z_k =
\begin{cases}
  1, & \text{if }\ \CI_c\geq \CI_k\\
  0, & \text{otherwise}.
\end{cases}
\]
This can be rewritten equivalently as
\[
z_k\in\{0,1\}\ \text{ , }\ z_k\Big[\sum_{q = 1}^{Q}(I_{q, c}- I_{q, k})w_{q}\Big]\geq 0
\]
and $z_{k}$ takes the maximal value of 0 and 1 such that the right inequality is fulfilled. Since the number $R_c$ of countries that are not better than country $c$ can now be written as $R_c = \sum_{k}z_k$ we obtain the optimization problem
\begin{align*}
  R_c^{*} &= \max\limits_{\vec{w},\vec{z}} \sum_{k = 1\atop k\neq c}^{C} z_k\\
  &\begin{aligned}
    \text{subject to} && z_k\Big[\sum\limits_{q = 1}^{Q}(I_{q, c}-I_{q, k})w_{q}\Big]& \geq 0,& k&=1,\dots,C,\ k\neq c\\
    && \sum\limits_{q=1}^Q w_q&= 10&&\\
    && w_q&\geq 0,& q&=1,\dots, Q\\
    && z_k & \in \{0,1\},& k&=1,\dots,C,\ k\neq c.
  \end{aligned}
\end{align*}
Note that this new formulation has the new auxiliary optimization variable $\vec{z}$ which is not of genuine interest but is only used to bring the problem to a simpler form. Indeed, that new objective function  is very simple in terms of the optimization variables; in fact it is a linear function. The only thing which makes this problem still difficult is the binary constraint $z_k\in\{0,1\}$.

To facilitate the problem further we use \emph{relaxation}, i.e. we relax the binary constraint $z_k\in\{0,1\}$ to the constraint $0\leq z_k\leq 1$. This enlarges the search space and could, in principle lead to a larger objective and a different solution. But note that the objective of the relaxed problem can not be larger since still $z_k\leq 1$ and also that the optimal solution of the relaxed problem
\begin{align*}
  &\max\limits_{\vec{w},\vec{z}} \sum_{k = 1\atop k\neq c}^{C} z_k\\
  &\begin{aligned}
    \st && z_k\Big[\sum\limits_{q = 1}^{Q}(I_{q, c}-I_{q, k})w_{q}\Big]& \geq 0,& k&=1,\dots,C,\ k\neq c\\
    && \sum\limits_{q=1}^Q w_q&= 10&&\\
    && w_q&\geq 0,& q&=1,\dots, Q\\
    && z_k & \in [0,1],& k&=1,\dots,C,\ k\neq c
  \end{aligned}
\end{align*} 
is indeed \emph{equal} to the solution of the original one. To see this, note that if $z_k > 0$ we obtain from the first constraint that $ \sum_{q = 1}^{Q}(I_{q, c}-I_{q, k})w_{q} \geq 0$ (i.e. country $c$ is not worse than country $k$). But since we maximize the objective function which is nothing else than the sum of the positive values of the $z_k$, we see that it is indeed optimal to set $z_k=1$ if $ \sum_{q = 1}^{Q}(I_{q, c}-I_{q, k})w_{q}\geq 0$.

With this formulation we have transformed the problem such that it is accessible to standard tools of mathematical optimization. 

\subsection{Second order rank maximization for country $c$ via distance}\label{subsec:secondorder}

The rank-optimal weights for country $c$ are in general not unique, i.e. there may be several (usually infinitely many) weightings which lead to the same maximal value $R_{c}^{*}$. However, the weightings that lead to the optimal rank are themselves of interest. In this section we describe a strategy that is to identify one weighting in the set of all rank-optimal weightings for country $c$ that obeys a certain characteristic. More precisely, we look for a rank-optimal weighting that maximizes the distance of the score of country $c$ to the score of the direct successor of country $c$ in an optimal ranking. In other words, we try to find the weights such that the country gets the highest possible rank in a way that it is most ahead of the next countries in the ranking.\footnote{\label{ftnt} Other second order criteria could also be considered, e.g. choosing the weights such that the rank is optimal, and the score is as large as possible, or the smallest weight is as large as possible, or the vector of weights is closest to a vector with equal weights. For the sake of brevity we do not consider other criteria here.}

We encounter this problem by introducing a new auxiliary variable $d$ representing this very distance and take this as new objective. This $d$ shall be non-negative (i.e. we constrain to $d\geq 0$) and $D_{c}^{*}$ shall be the largest possible number such that the score of country $c$ is still at least $D_{c}^{*}$ larger than the score of all countries which are dominated by country $c$. This is ensured by the constraints $z_{k}\left[\sum_{q = 1}^{Q}(I_{q, c} - I_{q, k})w_{q} - d\right] \geq 0$. To ensure that country $c$ still has the best possible rank, we reuse the known optimal rank $R_c^{*}$ of country $c$ and add the constraint, that the new weights should also lead to this rank. This can be done by constraining $\sum_{k = 1\atop k\neq c}^{C} z_{k} = R_{c}^{*}$. In total we obtain the following optimization problem:
\begin{align*}
  D_c^{*} &= \max\limits_{\vec{w},\vec{z}, d} d\\
  &\begin{aligned}
    \text{subject to} &&\sum_{k = 1\atop k\neq c}^{C} z_{k} &= R_{c}^{*}\\
    && z_k\Big[\sum\limits_{q = 1}^{Q}(I_{q, c}-I_{q, k})w_{q} - d\Big]& \geq 0,& k&=1,\dots,C,\ k\neq c\\
    && \sum\limits_{q=1}^Q w_q&= 10&&\\
    && w_q&\geq 0,& q&=1,\dots, Q\\
    && z_k & \in \{0,1\},& k&=1,\dots,C,\ k\neq c\\
    && d &\geq 0.
  \end{aligned}
\end{align*}

Any solution of this optimization problem produces weights $\vec{w}$
that lead to the optimal rank $R_c^{*}$ of the country $c$ in a way
that the score of country $c$ has the largest distance $D_c^{*}$ to the
countries that follow in the ranking. Note that the solution is in
general still not unique.

\subsection{First and second order rank maximization with integer weights}

In the OECD Better Life Index interactive tool, the user is asked to assign integer weights between 0 and 5 to the different dimensions of well-being. Therefore, it is natural to ask what the best possible rank for a country $c$ would be in case the weights $w_{q}$ were constrained to be integer numbers between 0 and 5. Indeed, we can answer this question if we adapt our model appropriately. For the purpose of first order rank maximization we simply change the domain of the weights $w_{q}$ accordingly. Moreover, we drop the previous constraint $\sum_{q=1}^Q w_q = 10$ and replace it by $\sum_{q=1}^Q w_q \geq 1$ which ensures again that at least one weight is non-zero. Apart from this, the first order optimization problem remains unchanged and the result is

\begin{align*}
  R_c^{*} &= \max\limits_{\vec{w},\vec{z}} \sum_{k = 1\atop k\neq c}^{C} z_k\\
  &\begin{aligned}
    \text{subject to} && z_k\Big[\sum\limits_{q = 1}^{Q}(I_{q, c}-I_{q, k})w_{q}\Big]& \geq 0,& k&=1,\dots,C,\ k\neq c\\
    && \sum\limits_{q=1}^Q w_q&\geq 1&&\\
    && 0\leq w_q&\leq 5,& q&=1,\dots, Q\\
    && w_q&\in\mathbb{Z},& q&=1,\dots, Q\\
    && z_k & \in \{0,1\},& k&=1,\dots,C,\ k\neq c.
  \end{aligned}
\end{align*}

In case of second order rank maximization we take into account that the interactive tool normalizes the corporate index for the country $c$ to a value $0\leq \CI_{c}\leq 10$. This is done by setting $\vec{\hat w} = 10 \cdot \vec{w} / \sum_{q}w_{q}$ and $\CI_c = \sum_{q=1}^Q \hat w_q I_{q,c}$. The distance between the scores of two countries is calculated accordingly. This aspect is implemented in the following model:

\begin{align*}
  D_c^{*} &= \max\limits_{\vec{w},\vec{z}, d} d\\
  &\begin{aligned}
    \text{subject to} &&\sum_{k = 1\atop k\neq c}^{C} z_{k} &= R_{c}^{*}\\
    && z_k\Big[\sum\limits_{q = 1}^{Q}(I_{q, c}-I_{q, k})\hat w_{q} - d\Big]& \geq 0,& k&=1,\dots,C,\ k\neq c\\
    && \sum\limits_{q=1}^Q w_q&\geq 1&&\\
    && 10 \cdot w_{q} \Big/ \Big(\sum\limits_{q=1}^Q w_q\Big) &= \hat w_{q},& q&=1,\dots, Q\\
    && 0\leq w_q&\leq 5,& q&=1,\dots, Q\\
    && w_q&\in\mathbb{Z},& q&=1,\dots, Q\\
    && z_k & \in \{0,1\},& k&=1,\dots,C,\ k\neq c\\
    && d &\geq 0.
  \end{aligned}
\end{align*}

Line by line the problem reads: We want to maximize the distance of country $c$ to its immediate followers in a ranking ($D_c^\ast$) over all weight vectors ($\vec{w}$), auxiliary variable vectors ($\vec{z}$), and auxiliary variable ($d$), subject to the constraints that the country still achieves its optimal rank, there is at least one positive weight, the weight vector can be appropriately normalized, weights are between zero and five, weights are integers, and auxiliary variables fulfill the constraints which they should fulfill by design.  

This last problem formulation is what we use in the following to compute and discuss rank-optimal weights. For comparison, the appendix includes results for the procedure without the restriction to integer weights from Subsection \ref{subsec:secondorder}. 

\subsection{Computation}

We used the freely available solver SCIP to solve all occuring optimization problems. As \citet{Achterberg2009SCIP} outlines, SCIP is a software framework for so called \textit{constraint integer programs} (CIPs). CIPs include problems with linear objective function, nonlinear constraints as well as both continuous and integer variables. An important feature of CIPs concerns the subproblems that arise after fixing all integer variables, in particular these problems must be \textit{linear programs}. The above formulated problems for first and second order rank maximization fit closely into this framework. The algorithmic framework of SCIP encompasses, among others, a \textit{branch-and-bound} approach, \textit{LP relaxation} and \textit{cutting plane separation}.

The dataset and MATLAB code which computes the results using SCIP is publicly available \citet{Lorenz.Brauer.ea2015RankoptimalWeighting}
It can be adjusted to compute rank-optimal weights for other composite indicators.
To reproduce the results in this paper, download and install the ``SCIP Optimization Suite'' from \url{http://scip.zib.de/} on your machine and download the data and m-files (MATLAB scripts) from \citet{Lorenz.Brauer.ea2015RankoptimalWeighting}.
The m-files that come with the data contain a driver file that imports the data, formats it for use with SCIP and calls SCIP to solve the optimization problems. Then, it presents the output as text format but also leaves the results in the workspace for further use. The whole calculation of the best ranks for all countries in OECD Better Life Index for 2013 (2014) took about 43 seconds (25 seconds) on a machine with  Intel\textsuperscript{\textregistered} Core\textsuperscript{\texttrademark} i7, 1.90Ghz processor and 3.7 GiB of RAM. Some subinstances took considerably longer to solve than others.

\section{Results}

Let us first show for some example countries how rank-optimal weights change the OECD Better Life Index in 2014. We compare these rankings with the ranking with equal weights. For the case of Germany, our procedure with integer weights outputs the following rank-optimal weights: 5 for Safety; 3 for Education; 2 for Income, Environment, and Work-Life Balance; 1 for Jobs; and 0 for all other dimensions (Housing, Community, Civic Engagement, Health, and Life Satisfaction). Using these weights, Table \ref{tab:DEPL} shows that Germany can jump from rank twelve, which it has with equal weights, to the top of the ranking slightly before Switzerland, Finland
and Denmark. If Germany wants to do a marketing campaign it could focus on people with preferences close to these weights because people with these preferences could feel in Germany to live in the ``best'' (or close to best) country for them. The distance to the second is quite small in this case (0.01 BLI points), therefore rank one is possibly only achievable with weights very close to these.  

This example also shows that sometimes a good mix of weights instead of a focus on just the countries best dimension enables a country to outperform all others and jump to the top. If Germany would for example focus the hypothetical marketing campaign solely on safety and education, it might come close to Poland's rank-optimal weights with 4 for Safety, 2 for Education and 0 for all other dimensions. The ranking for these weights is also shown in Table \ref{tab:DEPL}. Indeed, Poland jumps from rank 26 with equal weights to the top with these weights, while Germany is not even in the ``Top 5'' anymore. People with a focus first on Safety and second on Education (with safety having twice the importance of eduction) would thus have their ``best life'' in Poland. 

The example of Poland demonstrates that most countries achieve their best rank by setting the weights of many dimensions to zero. This somehow contradicts the idea of better life being multidimensional in essence. Let us point out that the computed rank-optimal weights need not necessarily be the only weights which bring the country to its best position. There might be other weight vectors which also produce this which have more (or even all) dimensions with non-zero weights. The weights our method produces are often to a large extent determined through our second order criterion which is to also achieve the largest possible advance to the immediate rival below in the ranking. Our method could be easily modified to incorporate minimal requirements to account for multidimensionality of a ``better life'' weight vector, but this would require a discussion and definition of these requirements. This is beyond the scope of this paper but we briefly discuss some possible modifications later.

\begin{table}
\caption{Top five places in the ranking with rank-optimal integer weights for \textbf{Germany} and \textbf{Poland} for the OECD Better Life Index data of 2014. Rank-optimal weights are shown on top of the ranking. Values in each dimension are shown as barplots. Color shadings of bars correspond to normalized weights for the dimension. Thus, white bars are weighted with zero for the Better Life Index. This is analog to the web-interface of the OECD.}
\label{tab:DEPL}
\centering
 \renewcommand{\arraystretch}{0}
\begin{tabular}{m{3mm} >{\bf}m{24mm} >{}m{5mm} >{\centering}m{6.2mm} m{57mm} }
 &  &  &  & \begin{tikzpicture}[scale=0.5] 
              \foreach \i/\y  in {1/Housing,2/Income,3/Jobs,4/Community,5/Education,6/Environment,7/Civic Engagement,8/Health,9/Life Satisfaction,10/Safety,11/Work Life Balance} {\draw (\i,0) node[anchor=south,rotate=90]{\y};}; 
              \end{tikzpicture}  \\ \midrule
                  & Rank-optimal weights \newline Germany 2014 & & & \begin{tikzpicture}[yscale=0.6/0.870967741935484,scale=0.5]
                                  \foreach \i/\y/\c/\w in {2/0.8/Cerulean/2,3/0.4/NavyBlue/1,5/1.2/LimeGreen/3,6/0.8/ForestGreen/2,10/2/Gray/5,11/0.8/Maroon/2} 
                                       {\draw[line width=10pt,\c] (\i,0)--(\i,\y)node[above,black,yshift=-6pt]{\w};};
                                  \draw (0.5,0)--(11.5,0);
                                  \end{tikzpicture}  \\ \midrule
 & Country & BLI & Rank eq.w. &  \\ \midrule
1 & Germany & 8.07 & 12 & \begin{tikzpicture}[yscale=1, scale=0.5]
                 \foreach \i/\y/\c/\f in {1/0.626/JungleGreen/0,2/0.527/Cerulean/40,3/0.826/NavyBlue/20,4/0.893/WildStrawberry/0,5/0.797/LimeGreen/60,6/0.877/ForestGreen/40,7/0.393/Dandelion/0,8/0.719/Plum/0,9/0.742/RedOrange/0,10/0.896/Gray/100,11/0.795/Maroon/40} 
                 {\draw[line width=0.5pt,black,fill=\c!\f] (\i-0.33,0)rectangle(\i+0.33,\y);};
                 \draw (0.5,0)--(11.5,0);
                 \end{tikzpicture} \\[3mm] 
2 & Switzerland & 8.06 & 4 & \begin{tikzpicture}[yscale=1, scale=0.5]
                 \foreach \i/\y/\c/\f in {1/0.623/JungleGreen/0,2/0.726/Cerulean/40,3/0.956/NavyBlue/20,4/0.929/WildStrawberry/0,5/0.745/LimeGreen/60,6/0.832/ForestGreen/40,7/0.337/Dandelion/0,8/0.925/Plum/0,9/1/RedOrange/0,10/0.87/Gray/100,11/0.713/Maroon/40} 
                 {\draw[line width=0.5pt,black,fill=\c!\f] (\i-0.33,0)rectangle(\i+0.33,\y);};
                 \draw (0.5,0)--(11.5,0);
                 \end{tikzpicture} \\[3mm] 
3 & Finland & 8.05 & 7 & \begin{tikzpicture}[yscale=1, scale=0.5]
                 \foreach \i/\y/\c/\f in {1/0.633/JungleGreen/0,2/0.349/Cerulean/40,3/0.748/NavyBlue/20,4/0.893/WildStrawberry/0,5/0.915/LimeGreen/60,6/0.9/ForestGreen/40,7/0.596/Dandelion/0,8/0.745/Plum/0,9/0.871/RedOrange/0,10/0.922/Gray/100,11/0.737/Maroon/40} 
                 {\draw[line width=0.5pt,black,fill=\c!\f] (\i-0.33,0)rectangle(\i+0.33,\y);};
                 \draw (0.5,0)--(11.5,0);
                 \end{tikzpicture} \\[3mm] 
4 & Denmark & 8.04 & 4 & \begin{tikzpicture}[yscale=1, scale=0.5]
                 \foreach \i/\y/\c/\f in {1/0.616/JungleGreen/0,2/0.396/Cerulean/40,3/0.81/NavyBlue/20,4/1/WildStrawberry/0,5/0.774/LimeGreen/60,6/0.9/ForestGreen/40,7/0.706/Dandelion/0,8/0.737/Plum/0,9/0.935/RedOrange/0,10/0.877/Gray/100,11/0.978/Maroon/40} 
                 {\draw[line width=0.5pt,black,fill=\c!\f] (\i-0.33,0)rectangle(\i+0.33,\y);};
                 \draw (0.5,0)--(11.5,0);
                 \end{tikzpicture} \\[3mm] 
5 & Canada & 8.02 & 4 & \begin{tikzpicture}[yscale=1, scale=0.5]
                 \foreach \i/\y/\c/\f in {1/0.766/JungleGreen/0,2/0.572/Cerulean/40,3/0.801/NavyBlue/20,4/0.929/WildStrawberry/0,5/0.765/LimeGreen/60,6/0.853/ForestGreen/40,7/0.584/Dandelion/0,8/0.918/Plum/0,9/0.935/RedOrange/0,10/0.972/Gray/100,11/0.613/Maroon/40} 
                 {\draw[line width=0.5pt,black,fill=\c!\f] (\i-0.33,0)rectangle(\i+0.33,\y);};
                 \draw (0.5,0)--(11.5,0);
                 \end{tikzpicture} \\[3mm] 
 & $\dots$ &  &  & $\dots$ \\ \midrule
\end{tabular}
 \renewcommand{\arraystretch}{0}
\begin{tabular}{m{3mm} >{\bf}m{24mm} >{}m{5mm} >{\centering}m{6.2mm} m{57mm} }
 &  &  &  & \begin{tikzpicture}[scale=0.5] 
              \foreach \i/\y  in {1/Housing,2/Income,3/Jobs,4/Community,5/Education,6/Environment,7/Civic Engagement,8/Health,9/Life Satisfaction,10/Safety,11/Work Life Balance} {\draw (\i,0) node[anchor=south,rotate=90]{\y};}; 
              \end{tikzpicture}  \\ \midrule
                  & Rank-optimal weights \newline Poland 2014 & & & \begin{tikzpicture}[yscale=0.6/0.870967741935484,scale=0.5]
                                  \foreach \i/\y/\c/\w in {5/1/LimeGreen/2,10/2/Gray/4} 
                                       {\draw[line width=10pt,\c] (\i,0)--(\i,\y)node[above,black,yshift=-6pt]{\w};};
                                  \draw (0.5,0)--(11.5,0);
                                  \end{tikzpicture}  \\ \midrule
 & Country & BLI & Rank eq.w. &  \\ \midrule
1 & Poland & 9.34 & 26 & \begin{tikzpicture}[yscale=1, scale=0.5]
                 \foreach \i/\y/\c/\f in {1/0.36/JungleGreen/0,2/0.129/Cerulean/0,3/0.524/NavyBlue/0,4/0.75/WildStrawberry/0,5/0.84/LimeGreen/50,6/0.487/ForestGreen/0,7/0.531/Dandelion/0,8/0.52/Plum/0,9/0.323/RedOrange/0,10/0.982/Gray/100,11/0.562/Maroon/0} 
                 {\draw[line width=0.5pt,black,fill=\c!\f] (\i-0.33,0)rectangle(\i+0.33,\y);};
                 \draw (0.5,0)--(11.5,0);
                 \end{tikzpicture} \\[3mm] 
2 & Japan & 9.24 & 20 & \begin{tikzpicture}[yscale=1, scale=0.5]
                 \foreach \i/\y/\c/\f in {1/0.484/JungleGreen/0,2/0.569/Cerulean/0,3/0.799/NavyBlue/0,4/0.786/WildStrawberry/0,5/0.782/LimeGreen/50,6/0.694/ForestGreen/0,7/0.393/Dandelion/0,8/0.496/Plum/0,9/0.419/RedOrange/0,10/0.996/Gray/100,11/0.526/Maroon/0} 
                 {\draw[line width=0.5pt,black,fill=\c!\f] (\i-0.33,0)rectangle(\i+0.33,\y);};
                 \draw (0.5,0)--(11.5,0);
                 \end{tikzpicture} \\[3mm] 
3 & Finland & 9.2 & 7 & \begin{tikzpicture}[yscale=1, scale=0.5]
                 \foreach \i/\y/\c/\f in {1/0.633/JungleGreen/0,2/0.349/Cerulean/0,3/0.748/NavyBlue/0,4/0.893/WildStrawberry/0,5/0.915/LimeGreen/50,6/0.9/ForestGreen/0,7/0.596/Dandelion/0,8/0.745/Plum/0,9/0.871/RedOrange/0,10/0.922/Gray/100,11/0.737/Maroon/0} 
                 {\draw[line width=0.5pt,black,fill=\c!\f] (\i-0.33,0)rectangle(\i+0.33,\y);};
                 \draw (0.5,0)--(11.5,0);
                 \end{tikzpicture} \\[3mm] 
4 & Canada & 9.03 & 4 & \begin{tikzpicture}[yscale=1, scale=0.5]
                 \foreach \i/\y/\c/\f in {1/0.766/JungleGreen/0,2/0.572/Cerulean/0,3/0.801/NavyBlue/0,4/0.929/WildStrawberry/0,5/0.765/LimeGreen/50,6/0.853/ForestGreen/0,7/0.584/Dandelion/0,8/0.918/Plum/0,9/0.935/RedOrange/0,10/0.972/Gray/100,11/0.613/Maroon/0} 
                 {\draw[line width=0.5pt,black,fill=\c!\f] (\i-0.33,0)rectangle(\i+0.33,\y);};
                 \draw (0.5,0)--(11.5,0);
                 \end{tikzpicture} \\[3mm] 
5 & Korea & 8.99 & 25 & \begin{tikzpicture}[yscale=1, scale=0.5]
                 \foreach \i/\y/\c/\f in {1/0.574/JungleGreen/0,2/0.229/Cerulean/0,3/0.759/NavyBlue/0,4/0.321/WildStrawberry/0,5/0.799/LimeGreen/50,6/0.537/ForestGreen/0,7/0.749/Dandelion/0,8/0.497/Plum/0,9/0.419/RedOrange/0,10/0.949/Gray/100,11/0.417/Maroon/0} 
                 {\draw[line width=0.5pt,black,fill=\c!\f] (\i-0.33,0)rectangle(\i+0.33,\y);};
                 \draw (0.5,0)--(11.5,0);
                 \end{tikzpicture} \\[3mm] 
 & $\dots$ &  &  & $\dots$ \\ \midrule
\end{tabular}
\end{table}

Table \ref{tab:ESAT} shows two other examples, Spain and Austria. Spain can jump twenty ranks upward to the top in an OECD world where better life is focusing on Health and slightly less on Work-life Balance (in proportion five to four). Austria jumps from rank 16 to rank 2 with weights 4 for Jobs and Community, 2 for Income and 1 for Life Satisfaction, but it can not overtake Switzerland with these weights. It could overtake Switzerland by putting all weight on Community, Civic Engagement and Safety, because it scores better than Switzerland in these three dimensions, but this would come with the cost of being overtaken by other countries. We could compute who would overtake, but we know already that there are no weights which could bring Austria to the top, because such weights would have come out in the optimization. 

Tables \ref{tab:2013d} and \ref{tab:2014d} show all discrete rank-optimal weights for the OECD Better Life Index in 2013 and 2014 respectively in a graphic form. The original numbers are available in  \citet{Lorenz.Brauer.ea2015RankoptimalWeighting}.
The table also shows for each country the rank which is achieved with these weights and the distance in Better Life Index to the next worse in this ranking. The maximization of the latter is the second order criterion in the optimization problem. Countries are sorted according to these two values. It is visible that more than half of the 36 countries (19 in 2014, 17 in 2013) can be brought to the (unshared) top of the ranking with appropriate weights. All these countries can claim to have the ``best'' life given a particular definition of a better life in terms of weights in the framework of the OECD Better Life Index (which also allows to define ``better life'' based on a single non-zero dimension). All countries can make substantial improvements in ranks with respect
to the equal weighting through their rank optimal weights. Australia, which is on top with equal weights, is obviously not able to improve the rank, but it improves its distance to the second one through rank-optimal weighting. Nevertheless, there are differences in improvement. For example Mexico can improve from the second lowest rank 35 to rank 7 in 2014 (rank 10 in 2013), while Chile can only improve from rank 34 to 21. The high potential of Mexico to jump comes through the fact that Mexico is very weak in most dimensions, but fairly good in Life Satisfaction. In comparison, Chile is relatively weak in all dimensions, which limits the possibility to overtake other countries.  

\begin{table}
\caption{Top five places in the ranking with rank-optimal integer weights for \textbf{Spain} and \textbf{Austria} for the OECD Better Life Index data of 2014. Rank-optimal weights are shown on top of the ranking. Values in each dimension are shown as barplots. Color shadings of bars correspond to normalized weights for the dimension. Thus, white bars are weighted with zero for the Better Life Index. This is analog to the web-interface of the OECD.}
\label{tab:ESAT}
\centering
 \renewcommand{\arraystretch}{0}
\begin{tabular}{m{3mm} >{\bf}m{24mm} >{}m{5mm} >{\centering}m{6.2mm} m{57mm} }
 &  &  &  & \begin{tikzpicture}[scale=0.5] 
              \foreach \i/\y  in {1/Housing,2/Income,3/Jobs,4/Community,5/Education,6/Environment,7/Civic Engagement,8/Health,9/Life Satisfaction,10/Safety,11/Work Life Balance} {\draw (\i,0) node[anchor=south,rotate=90]{\y};}; 
              \end{tikzpicture}  \\ \midrule
                  & Rank-optimal weights \newline Spain 2014 & & & \begin{tikzpicture}[yscale=0.6/0.779606625258799,scale=0.5]
                                  \foreach \i/\y/\c/\w in {8/2/Plum/5,11/1.6/Maroon/4} 
                                       {\draw[line width=10pt,\c] (\i,0)--(\i,\y)node[above,black,yshift=-6pt]{\w};};
                                  \draw (0.5,0)--(11.5,0);
                                  \end{tikzpicture}  \\ \midrule
 & Country & BLI & Rank eq.w. &  \\ \midrule
1 & Spain & 8.93 & 21 & \begin{tikzpicture}[yscale=1, scale=0.5]
                 \foreach \i/\y/\c/\f in {1/0.688/JungleGreen/0,2/0.293/Cerulean/0,3/0.258/NavyBlue/0,4/0.857/WildStrawberry/0,5/0.536/LimeGreen/0,6/0.59/ForestGreen/0,7/0.506/Dandelion/0,8/0.861/Plum/100,9/0.484/RedOrange/0,10/0.866/Gray/0,11/0.933/Maroon/80} 
                 {\draw[line width=0.5pt,black,fill=\c!\f] (\i-0.33,0)rectangle(\i+0.33,\y);};
                 \draw (0.5,0)--(11.5,0);
                 \end{tikzpicture} \\[3mm] 
2 & Netherlands & 8.51 & 7 & \begin{tikzpicture}[yscale=1, scale=0.5]
                 \foreach \i/\y/\c/\f in {1/0.688/JungleGreen/0,2/0.525/Cerulean/0,3/0.868/NavyBlue/0,4/0.857/WildStrawberry/0,5/0.761/LimeGreen/0,6/0.688/ForestGreen/0,7/0.511/Dandelion/0,8/0.829/Plum/100,9/0.871/RedOrange/0,10/0.832/Gray/0,11/0.878/Maroon/80} 
                 {\draw[line width=0.5pt,black,fill=\c!\f] (\i-0.33,0)rectangle(\i+0.33,\y);};
                 \draw (0.5,0)--(11.5,0);
                 \end{tikzpicture} \\[3mm] 
3 & Sweden & 8.51 & 2 & \begin{tikzpicture}[yscale=1, scale=0.5]
                 \foreach \i/\y/\c/\f in {1/0.625/JungleGreen/0,2/0.496/Cerulean/0,3/0.782/NavyBlue/0,4/0.821/WildStrawberry/0,5/0.79/LimeGreen/0,6/0.986/ForestGreen/0,7/0.878/Dandelion/0,8/0.884/Plum/100,9/0.871/RedOrange/0,10/0.821/Gray/0,11/0.809/Maroon/80} 
                 {\draw[line width=0.5pt,black,fill=\c!\f] (\i-0.33,0)rectangle(\i+0.33,\y);};
                 \draw (0.5,0)--(11.5,0);
                 \end{tikzpicture} \\[3mm] 
4 & Denmark & 8.44 & 4 & \begin{tikzpicture}[yscale=1, scale=0.5]
                 \foreach \i/\y/\c/\f in {1/0.616/JungleGreen/0,2/0.396/Cerulean/0,3/0.81/NavyBlue/0,4/1/WildStrawberry/0,5/0.774/LimeGreen/0,6/0.9/ForestGreen/0,7/0.706/Dandelion/0,8/0.737/Plum/100,9/0.935/RedOrange/0,10/0.877/Gray/0,11/0.978/Maroon/80} 
                 {\draw[line width=0.5pt,black,fill=\c!\f] (\i-0.33,0)rectangle(\i+0.33,\y);};
                 \draw (0.5,0)--(11.5,0);
                 \end{tikzpicture} \\[3mm] 
5 & Norway & 8.36 & 2 & \begin{tikzpicture}[yscale=1, scale=0.5]
                 \foreach \i/\y/\c/\f in {1/0.764/JungleGreen/0,2/0.392/Cerulean/0,3/0.92/NavyBlue/0,4/0.893/WildStrawberry/0,5/0.719/LimeGreen/0,6/0.896/ForestGreen/0,7/0.651/Dandelion/0,8/0.808/Plum/100,9/0.968/RedOrange/0,10/0.873/Gray/0,11/0.871/Maroon/80} 
                 {\draw[line width=0.5pt,black,fill=\c!\f] (\i-0.33,0)rectangle(\i+0.33,\y);};
                 \draw (0.5,0)--(11.5,0);
                 \end{tikzpicture} \\[3mm] 
 & $\dots$ &  &  & $\dots$ \\ \midrule
\end{tabular}
 \renewcommand{\arraystretch}{0}
\begin{tabular}{m{3mm} >{\bf}m{24mm} >{}m{5mm} >{\centering}m{6.2mm} m{57mm} }
 &  &  &  & \begin{tikzpicture}[scale=0.5] 
              \foreach \i/\y  in {1/Housing,2/Income,3/Jobs,4/Community,5/Education,6/Environment,7/Civic Engagement,8/Health,9/Life Satisfaction,10/Safety,11/Work Life Balance} {\draw (\i,0) node[anchor=south,rotate=90]{\y};}; 
              \end{tikzpicture}  \\ \midrule
                  & Rank-optimal weights \newline Austria 2014 & & & \begin{tikzpicture}[yscale=0.6/0.873826777087647,scale=0.5]
                                  \foreach \i/\y/\c/\w in {2/1/Cerulean/2,3/2/NavyBlue/4,4/2/WildStrawberry/4,9/0.5/RedOrange/1} 
                                       {\draw[line width=10pt,\c] (\i,0)--(\i,\y)node[above,black,yshift=-6pt]{\w};};
                                  \draw (0.5,0)--(11.5,0);
                                  \end{tikzpicture}  \\ \midrule
 & Country & BLI & Rank eq.w. &  \\ \midrule
1 & Switzerland & 9.08 & 4 & \begin{tikzpicture}[yscale=1, scale=0.5]
                 \foreach \i/\y/\c/\f in {1/0.623/JungleGreen/0,2/0.726/Cerulean/50,3/0.956/NavyBlue/100,4/0.929/WildStrawberry/100,5/0.745/LimeGreen/0,6/0.832/ForestGreen/0,7/0.337/Dandelion/0,8/0.925/Plum/0,9/1/RedOrange/25,10/0.87/Gray/0,11/0.713/Maroon/0} 
                 {\draw[line width=0.5pt,black,fill=\c!\f] (\i-0.33,0)rectangle(\i+0.33,\y);};
                 \draw (0.5,0)--(11.5,0);
                 \end{tikzpicture} \\[3mm] 
2 & Austria & 8.36 & 15 & \begin{tikzpicture}[yscale=1, scale=0.5]
                 \foreach \i/\y/\c/\f in {1/0.582/JungleGreen/0,2/0.497/Cerulean/50,3/0.859/NavyBlue/100,4/0.964/WildStrawberry/100,5/0.669/LimeGreen/0,6/0.738/ForestGreen/0,7/0.564/Dandelion/0,8/0.763/Plum/0,9/0.903/RedOrange/25,10/0.905/Gray/0,11/0.599/Maroon/0} 
                 {\draw[line width=0.5pt,black,fill=\c!\f] (\i-0.33,0)rectangle(\i+0.33,\y);};
                 \draw (0.5,0)--(11.5,0);
                 \end{tikzpicture} \\[3mm] 
3 & Iceland & 8.26 & 10 & \begin{tikzpicture}[yscale=1, scale=0.5]
                 \foreach \i/\y/\c/\f in {1/0.595/JungleGreen/0,2/0.361/Cerulean/50,3/0.866/NavyBlue/100,4/1/WildStrawberry/100,5/0.726/LimeGreen/0,6/0.878/ForestGreen/0,7/0.527/Dandelion/0,8/0.886/Plum/0,9/0.903/RedOrange/25,10/0.919/Gray/0,11/0.568/Maroon/0} 
                 {\draw[line width=0.5pt,black,fill=\c!\f] (\i-0.33,0)rectangle(\i+0.33,\y);};
                 \draw (0.5,0)--(11.5,0);
                 \end{tikzpicture} \\[3mm] 
4 & United States & 8.24 & 7 & \begin{tikzpicture}[yscale=1, scale=0.5]
                 \foreach \i/\y/\c/\f in {1/0.789/JungleGreen/0,2/1/Cerulean/50,3/0.796/NavyBlue/100,4/0.786/WildStrawberry/100,5/0.699/LimeGreen/0,6/0.784/ForestGreen/0,7/0.536/Dandelion/0,8/0.851/Plum/0,9/0.742/RedOrange/25,10/0.894/Gray/0,11/0.53/Maroon/0} 
                 {\draw[line width=0.5pt,black,fill=\c!\f] (\i-0.33,0)rectangle(\i+0.33,\y);};
                 \draw (0.5,0)--(11.5,0);
                 \end{tikzpicture} \\[3mm] 
5 & Norway & 8.19 & 2 & \begin{tikzpicture}[yscale=1, scale=0.5]
                 \foreach \i/\y/\c/\f in {1/0.764/JungleGreen/0,2/0.392/Cerulean/50,3/0.92/NavyBlue/100,4/0.893/WildStrawberry/100,5/0.719/LimeGreen/0,6/0.896/ForestGreen/0,7/0.651/Dandelion/0,8/0.808/Plum/0,9/0.968/RedOrange/25,10/0.873/Gray/0,11/0.871/Maroon/0} 
                 {\draw[line width=0.5pt,black,fill=\c!\f] (\i-0.33,0)rectangle(\i+0.33,\y);};
                 \draw (0.5,0)--(11.5,0);
                 \end{tikzpicture} \\[3mm] 
 & $\dots$ &  &  & $\dots$ \\ \midrule
\end{tabular}
\end{table}

\begin{table}
\centering
\caption{Rank-optimal integer weights (0,1,2,3,4,5) for the OECD Better Life Index 2013.}
\label{tab:2013d}
\renewcommand{\arraystretch}{0}
\begin{tabular}{m{30mm} >{\centering\bf}m{3mm} >{\centering}m{5mm} >{\centering}m{3mm} m{57mm} }
\hline\noalign{\smallskip}
Country & \rotatebox{90}{Top rank} & \rotatebox{90}{Distance to next} & \rotatebox{90}{Rank equal weights} &\begin{tikzpicture}[scale=0.5] 
              \foreach \i/\y  in {1/Housing,2/Income,3/Jobs,4/Community,5/Education,6/Environment,7/Civic Engagement,8/Health,9/Life Satisfaction,10/Safety,11/Work Life Balance} {\draw (\i,0) node[anchor=south,rotate=90]{\y};}; 
              \end{tikzpicture}  \\ \midrule
United States & 1 & 2.189 & 6 & \begin{tikzpicture}[yscale=0.6/1,scale=0.5]
                                  \foreach \i/\y/\c in {2/1/Cerulean} {\draw[line width=10pt,\c] (\i,0)--(\i,\y);};
                                  \draw (0.5,0)--(11.5,0);
                                  \end{tikzpicture} \\[-2.2mm] 
Switzerland & 1 & 1.227 & 5 & \begin{tikzpicture}[yscale=0.6/3,scale=0.5]
                                  \foreach \i/\y/\c in {2/1/Cerulean,9/3/RedOrange} {\draw[line width=10pt,\c] (\i,0)--(\i,\y);};
                                  \draw (0.5,0)--(11.5,0);
                                  \end{tikzpicture} \\[-2.2mm] 
Australia & 1 & 1.055 & 1 & \begin{tikzpicture}[yscale=0.6/5,scale=0.5]
                                  \foreach \i/\y/\c in {1/5/JungleGreen,7/5/Dandelion,10/4/Gray} {\draw[line width=10pt,\c] (\i,0)--(\i,\y);};
                                  \draw (0.5,0)--(11.5,0);
                                  \end{tikzpicture} \\[-2.2mm] 
Iceland & 1 & 1.039 & 9 & \begin{tikzpicture}[yscale=0.6/3,scale=0.5]
                                  \foreach \i/\y/\c in {4/3/WildStrawberry,9/1/RedOrange} {\draw[line width=10pt,\c] (\i,0)--(\i,\y);};
                                  \draw (0.5,0)--(11.5,0);
                                  \end{tikzpicture} \\[-2.2mm] 
Finland & 1 & 0.883 & 12 & \begin{tikzpicture}[yscale=0.6/5,scale=0.5]
                                  \foreach \i/\y/\c in {5/5/LimeGreen,11/1/Maroon} {\draw[line width=10pt,\c] (\i,0)--(\i,\y);};
                                  \draw (0.5,0)--(11.5,0);
                                  \end{tikzpicture} \\[-2.2mm] 
Sweden & 1 & 0.736 & 1 & \begin{tikzpicture}[yscale=0.6/4,scale=0.5]
                                  \foreach \i/\y/\c in {5/3/LimeGreen,6/4/ForestGreen,7/3/Dandelion,9/4/RedOrange,11/1/Maroon} {\draw[line width=10pt,\c] (\i,0)--(\i,\y);};
                                  \draw (0.5,0)--(11.5,0);
                                  \end{tikzpicture} \\[-2.2mm] 
Norway & 1 & 0.668 & 3 & \begin{tikzpicture}[yscale=0.6/4,scale=0.5]
                                  \foreach \i/\y/\c in {1/3/JungleGreen,3/4/NavyBlue,6/3/ForestGreen,11/2/Maroon} {\draw[line width=10pt,\c] (\i,0)--(\i,\y);};
                                  \draw (0.5,0)--(11.5,0);
                                  \end{tikzpicture} \\[-2.2mm] 
Denmark & 1 & 0.559 & 7 & \begin{tikzpicture}[yscale=0.6/5,scale=0.5]
                                  \foreach \i/\y/\c in {7/1/Dandelion,11/5/Maroon} {\draw[line width=10pt,\c] (\i,0)--(\i,\y);};
                                  \draw (0.5,0)--(11.5,0);
                                  \end{tikzpicture} \\[-2.2mm] 
Japan & 1 & 0.525 & 21 & \begin{tikzpicture}[yscale=0.6/5,scale=0.5]
                                  \foreach \i/\y/\c in {2/1/Cerulean,5/3/LimeGreen,10/5/Gray} {\draw[line width=10pt,\c] (\i,0)--(\i,\y);};
                                  \draw (0.5,0)--(11.5,0);
                                  \end{tikzpicture} \\[-2.2mm] 
United Kingdom & 1 & 0.479 & 9 & \begin{tikzpicture}[yscale=0.6/4,scale=0.5]
                                  \foreach \i/\y/\c in {2/1/Cerulean,6/4/ForestGreen,10/3/Gray} {\draw[line width=10pt,\c] (\i,0)--(\i,\y);};
                                  \draw (0.5,0)--(11.5,0);
                                  \end{tikzpicture} \\[-2.2mm] 
Ireland & 1 & 0.479 & 14 & \begin{tikzpicture}[yscale=0.6/4,scale=0.5]
                                  \foreach \i/\y/\c in {1/3/JungleGreen,4/4/WildStrawberry,11/1/Maroon} {\draw[line width=10pt,\c] (\i,0)--(\i,\y);};
                                  \draw (0.5,0)--(11.5,0);
                                  \end{tikzpicture} \\[-2.2mm] 
Luxembourg & 1 & 0.449 & 12 & \begin{tikzpicture}[yscale=0.6/5,scale=0.5]
                                  \foreach \i/\y/\c in {2/3/Cerulean,3/5/NavyBlue,7/3/Dandelion,11/4/Maroon} {\draw[line width=10pt,\c] (\i,0)--(\i,\y);};
                                  \draw (0.5,0)--(11.5,0);
                                  \end{tikzpicture} \\[-2.2mm] 
Netherlands & 1 & 0.349 & 7 & \begin{tikzpicture}[yscale=0.6/5,scale=0.5]
                                  \foreach \i/\y/\c in {1/3/JungleGreen,2/2/Cerulean,3/1/NavyBlue,4/5/WildStrawberry,11/5/Maroon} {\draw[line width=10pt,\c] (\i,0)--(\i,\y);};
                                  \draw (0.5,0)--(11.5,0);
                                  \end{tikzpicture} \\[-2.2mm] 
Canada & 1 & 0.344 & 3 & \begin{tikzpicture}[yscale=0.6/4,scale=0.5]
                                  \foreach \i/\y/\c in {1/3/JungleGreen,2/1/Cerulean,4/2/WildStrawberry,9/2/RedOrange,10/4/Gray} {\draw[line width=10pt,\c] (\i,0)--(\i,\y);};
                                  \draw (0.5,0)--(11.5,0);
                                  \end{tikzpicture} \\[-2.2mm] 
Belgium & 1 & 0.1 & 16 & \begin{tikzpicture}[yscale=0.6/5,scale=0.5]
                                  \foreach \i/\y/\c in {2/3/Cerulean,5/2/LimeGreen,7/1/Dandelion,11/5/Maroon} {\draw[line width=10pt,\c] (\i,0)--(\i,\y);};
                                  \draw (0.5,0)--(11.5,0);
                                  \end{tikzpicture} \\[-2.2mm] 
New Zealand & 1 & 0.081 & 9 & \begin{tikzpicture}[yscale=0.6/5,scale=0.5]
                                  \foreach \i/\y/\c in {6/1/ForestGreen,8/5/Plum,10/1/Gray,11/1/Maroon} {\draw[line width=10pt,\c] (\i,0)--(\i,\y);};
                                  \draw (0.5,0)--(11.5,0);
                                  \end{tikzpicture} \\[-2.2mm] 
Spain & 1 & 0.075 & 20 & \begin{tikzpicture}[yscale=0.6/5,scale=0.5]
                                  \foreach \i/\y/\c in {8/5/Plum,10/2/Gray,11/4/Maroon} {\draw[line width=10pt,\c] (\i,0)--(\i,\y);};
                                  \draw (0.5,0)--(11.5,0);
                                  \end{tikzpicture} \\[-2.2mm] 
Poland & 2 & 0.019 & 25 & \begin{tikzpicture}[yscale=0.6/5,scale=0.5]
                                  \foreach \i/\y/\c in {10/5/Gray,11/1/Maroon} {\draw[line width=10pt,\c] (\i,0)--(\i,\y);};
                                  \draw (0.5,0)--(11.5,0);
                                  \end{tikzpicture} \\[-2.2mm] 
Korea & 2 & 0.006 & 26 & \begin{tikzpicture}[yscale=0.6/5,scale=0.5]
                                  \foreach \i/\y/\c in {5/2/LimeGreen,7/3/Dandelion,10/5/Gray} {\draw[line width=10pt,\c] (\i,0)--(\i,\y);};
                                  \draw (0.5,0)--(11.5,0);
                                  \end{tikzpicture} \\[-2.2mm] 
Germany & 3 & 0.076 & 17 & \begin{tikzpicture}[yscale=0.6/5,scale=0.5]
                                  \foreach \i/\y/\c in {2/1/Cerulean,5/5/LimeGreen,6/5/ForestGreen} {\draw[line width=10pt,\c] (\i,0)--(\i,\y);};
                                  \draw (0.5,0)--(11.5,0);
                                  \end{tikzpicture} \\[-2.2mm] 
Austria & 3 & 0.076 & 14 & \begin{tikzpicture}[yscale=0.6/5,scale=0.5]
                                  \foreach \i/\y/\c in {3/5/NavyBlue,7/2/Dandelion,9/1/RedOrange,10/3/Gray} {\draw[line width=10pt,\c] (\i,0)--(\i,\y);};
                                  \draw (0.5,0)--(11.5,0);
                                  \end{tikzpicture} \\[-2.2mm] 
Israel & 5 & 0.128 & 24 & \begin{tikzpicture}[yscale=0.6/1,scale=0.5]
                                  \foreach \i/\y/\c in {8/1/Plum} {\draw[line width=10pt,\c] (\i,0)--(\i,\y);};
                                  \draw (0.5,0)--(11.5,0);
                                  \end{tikzpicture} \\[-2.2mm] 
Slovenia & 6 & 0.065 & 19 & \begin{tikzpicture}[yscale=0.6/5,scale=0.5]
                                  \foreach \i/\y/\c in {4/1/WildStrawberry,5/5/LimeGreen,7/5/Dandelion,11/1/Maroon} {\draw[line width=10pt,\c] (\i,0)--(\i,\y);};
                                  \draw (0.5,0)--(11.5,0);
                                  \end{tikzpicture} \\[-2.2mm] 
Russian Federation & 7 & 0.019 & 32 & \begin{tikzpicture}[yscale=0.6/1,scale=0.5]
                                  \foreach \i/\y/\c in {11/1/Maroon} {\draw[line width=10pt,\c] (\i,0)--(\i,\y);};
                                  \draw (0.5,0)--(11.5,0);
                                  \end{tikzpicture} \\[-2.2mm] 
France & 8 & 0.079 & 18 & \begin{tikzpicture}[yscale=0.6/4,scale=0.5]
                                  \foreach \i/\y/\c in {1/1/JungleGreen,2/4/Cerulean,4/1/WildStrawberry,8/1/Plum,11/1/Maroon} {\draw[line width=10pt,\c] (\i,0)--(\i,\y);};
                                  \draw (0.5,0)--(11.5,0);
                                  \end{tikzpicture} \\[-2.2mm] 
Hungary & 8 & 0.006 & 29 & \begin{tikzpicture}[yscale=0.6/5,scale=0.5]
                                  \foreach \i/\y/\c in {10/4/Gray,11/5/Maroon} {\draw[line width=10pt,\c] (\i,0)--(\i,\y);};
                                  \draw (0.5,0)--(11.5,0);
                                  \end{tikzpicture} \\[-2.2mm] 
Portugal & 9 & 0.24 & 28 & \begin{tikzpicture}[yscale=0.6/5,scale=0.5]
                                  \foreach \i/\y/\c in {1/5/JungleGreen} {\draw[line width=10pt,\c] (\i,0)--(\i,\y);};
                                  \draw (0.5,0)--(11.5,0);
                                  \end{tikzpicture} \\[-2.2mm] 
Estonia & 9 & 0.033 & 31 & \begin{tikzpicture}[yscale=0.6/1,scale=0.5]
                                  \foreach \i/\y/\c in {5/1/LimeGreen} {\draw[line width=10pt,\c] (\i,0)--(\i,\y);};
                                  \draw (0.5,0)--(11.5,0);
                                  \end{tikzpicture} \\[-2.2mm] 
Mexico & 10 & 0.323 & 35 & \begin{tikzpicture}[yscale=0.6/1,scale=0.5]
                                  \foreach \i/\y/\c in {9/1/RedOrange} {\draw[line width=10pt,\c] (\i,0)--(\i,\y);};
                                  \draw (0.5,0)--(11.5,0);
                                  \end{tikzpicture} \\[-2.2mm] 
Italy & 10 & 0.075 & 22 & \begin{tikzpicture}[yscale=0.6/4,scale=0.5]
                                  \foreach \i/\y/\c in {7/1/Dandelion,8/1/Plum,11/4/Maroon} {\draw[line width=10pt,\c] (\i,0)--(\i,\y);};
                                  \draw (0.5,0)--(11.5,0);
                                  \end{tikzpicture} \\[-2.2mm] 
Czech Republic & 10 & 0.043 & 22 & \begin{tikzpicture}[yscale=0.6/5,scale=0.5]
                                  \foreach \i/\y/\c in {5/3/LimeGreen,10/5/Gray} {\draw[line width=10pt,\c] (\i,0)--(\i,\y);};
                                  \draw (0.5,0)--(11.5,0);
                                  \end{tikzpicture} \\[-2.2mm] 
Slovak Republic & 10 & 0.009 & 26 & \begin{tikzpicture}[yscale=0.6/5,scale=0.5]
                                  \foreach \i/\y/\c in {6/1/ForestGreen,10/5/Gray,11/2/Maroon} {\draw[line width=10pt,\c] (\i,0)--(\i,\y);};
                                  \draw (0.5,0)--(11.5,0);
                                  \end{tikzpicture} \\[-2.2mm] 
Turkey & 11 & 0.114 & 36 & \begin{tikzpicture}[yscale=0.6/1,scale=0.5]
                                  \foreach \i/\y/\c in {7/1/Dandelion} {\draw[line width=10pt,\c] (\i,0)--(\i,\y);};
                                  \draw (0.5,0)--(11.5,0);
                                  \end{tikzpicture} \\[-2.2mm] 
Greece & 12 & 0.036 & 29 & \begin{tikzpicture}[yscale=0.6/4,scale=0.5]
                                  \foreach \i/\y/\c in {8/4/Plum,10/3/Gray} {\draw[line width=10pt,\c] (\i,0)--(\i,\y);};
                                  \draw (0.5,0)--(11.5,0);
                                  \end{tikzpicture} \\[-2.2mm] 
Brazil & 18 & 0.094 & 33 & \begin{tikzpicture}[yscale=0.6/3,scale=0.5]
                                  \foreach \i/\y/\c in {7/2/Dandelion,9/3/RedOrange} {\draw[line width=10pt,\c] (\i,0)--(\i,\y);};
                                  \draw (0.5,0)--(11.5,0);
                                  \end{tikzpicture} \\[-2.2mm] 
Chile & 21 & 0.05 & 34 & \begin{tikzpicture}[yscale=0.6/5,scale=0.5]
                                  \foreach \i/\y/\c in {9/5/RedOrange,10/1/Gray} {\draw[line width=10pt,\c] (\i,0)--(\i,\y);};
                                  \draw (0.5,0)--(11.5,0);
                                  \end{tikzpicture} \\[-2.2mm] 
\noalign{\smallskip}\hline
\end{tabular}
\end{table}

\begin{table}
\centering
\caption{Rank-optimal integer weights (0,1,2,3,4,5) for the OECD Better Life Index 2014.}
\label{tab:2014d}
\renewcommand{\arraystretch}{0}
\begin{tabular}{m{30mm} >{\centering\bf}m{3mm} >{\centering}m{5mm} >{\centering}m{3mm} m{57mm} }
\hline\noalign{\smallskip}
Country & \rotatebox{90}{Top rank} & \rotatebox{90}{Distance to next} & \rotatebox{90}{Rank equal weights} &\begin{tikzpicture}[scale=0.5] 
              \foreach \i/\y  in {1/Housing,2/Income,3/Jobs,4/Community,5/Education,6/Environment,7/Civic Engagement,8/Health,9/Life Satisfaction,10/Safety,11/Work Life Balance} {\draw (\i,0) node[anchor=south,rotate=90]{\y};}; 
              \end{tikzpicture}  \\ \midrule
United States & 1 & 2.739 & 7 & \begin{tikzpicture}[yscale=0.6/1,scale=0.5]
                                  \foreach \i/\y/\c in {2/1/Cerulean} {\draw[line width=10pt,\c] (\i,0)--(\i,\y);};
                                  \draw (0.5,0)--(11.5,0);
                                  \end{tikzpicture} \\[-2.2mm] 
Australia & 1 & 1.086 & 1 & \begin{tikzpicture}[yscale=0.6/2,scale=0.5]
                                  \foreach \i/\y/\c in {1/2/JungleGreen,7/2/Dandelion,10/2/Gray} {\draw[line width=10pt,\c] (\i,0)--(\i,\y);};
                                  \draw (0.5,0)--(11.5,0);
                                  \end{tikzpicture} \\[-2.2mm] 
Finland & 1 & 1.046 & 7 & \begin{tikzpicture}[yscale=0.6/5,scale=0.5]
                                  \foreach \i/\y/\c in {5/5/LimeGreen,9/1/RedOrange} {\draw[line width=10pt,\c] (\i,0)--(\i,\y);};
                                  \draw (0.5,0)--(11.5,0);
                                  \end{tikzpicture} \\[-2.2mm] 
Denmark & 1 & 1.023 & 4 & \begin{tikzpicture}[yscale=0.6/5,scale=0.5]
                                  \foreach \i/\y/\c in {4/5/WildStrawberry,7/1/Dandelion,11/5/Maroon} {\draw[line width=10pt,\c] (\i,0)--(\i,\y);};
                                  \draw (0.5,0)--(11.5,0);
                                  \end{tikzpicture} \\[-2.2mm] 
Switzerland & 1 & 1.021 & 4 & \begin{tikzpicture}[yscale=0.6/4,scale=0.5]
                                  \foreach \i/\y/\c in {2/2/Cerulean,3/3/NavyBlue,8/1/Plum,9/4/RedOrange} {\draw[line width=10pt,\c] (\i,0)--(\i,\y);};
                                  \draw (0.5,0)--(11.5,0);
                                  \end{tikzpicture} \\[-2.2mm] 
Sweden & 1 & 0.783 & 2 & \begin{tikzpicture}[yscale=0.6/4,scale=0.5]
                                  \foreach \i/\y/\c in {6/4/ForestGreen} {\draw[line width=10pt,\c] (\i,0)--(\i,\y);};
                                  \draw (0.5,0)--(11.5,0);
                                  \end{tikzpicture} \\[-2.2mm] 
Norway & 1 & 0.757 & 2 & \begin{tikzpicture}[yscale=0.6/3,scale=0.5]
                                  \foreach \i/\y/\c in {1/3/JungleGreen,3/2/NavyBlue,6/1/ForestGreen,11/1/Maroon} {\draw[line width=10pt,\c] (\i,0)--(\i,\y);};
                                  \draw (0.5,0)--(11.5,0);
                                  \end{tikzpicture} \\[-2.2mm] 
Spain & 1 & 0.421 & 21 & \begin{tikzpicture}[yscale=0.6/5,scale=0.5]
                                  \foreach \i/\y/\c in {8/5/Plum,11/4/Maroon} {\draw[line width=10pt,\c] (\i,0)--(\i,\y);};
                                  \draw (0.5,0)--(11.5,0);
                                  \end{tikzpicture} \\[-2.2mm] 
New Zealand & 1 & 0.378 & 10 & \begin{tikzpicture}[yscale=0.6/5,scale=0.5]
                                  \foreach \i/\y/\c in {4/5/WildStrawberry,7/1/Dandelion,8/4/Plum,11/1/Maroon} {\draw[line width=10pt,\c] (\i,0)--(\i,\y);};
                                  \draw (0.5,0)--(11.5,0);
                                  \end{tikzpicture} \\[-2.2mm] 
Canada & 1 & 0.37 & 4 & \begin{tikzpicture}[yscale=0.6/5,scale=0.5]
                                  \foreach \i/\y/\c in {1/3/JungleGreen,2/1/Cerulean,9/4/RedOrange,10/5/Gray} {\draw[line width=10pt,\c] (\i,0)--(\i,\y);};
                                  \draw (0.5,0)--(11.5,0);
                                  \end{tikzpicture} \\[-2.2mm] 
Belgium & 1 & 0.319 & 12 & \begin{tikzpicture}[yscale=0.6/5,scale=0.5]
                                  \foreach \i/\y/\c in {1/3/JungleGreen,2/3/Cerulean,11/5/Maroon} {\draw[line width=10pt,\c] (\i,0)--(\i,\y);};
                                  \draw (0.5,0)--(11.5,0);
                                  \end{tikzpicture} \\[-2.2mm] 
Iceland & 1 & 0.232 & 10 & \begin{tikzpicture}[yscale=0.6/2,scale=0.5]
                                  \foreach \i/\y/\c in {3/1/NavyBlue,4/2/WildStrawberry,10/1/Gray} {\draw[line width=10pt,\c] (\i,0)--(\i,\y);};
                                  \draw (0.5,0)--(11.5,0);
                                  \end{tikzpicture} \\[-2.2mm] 
Japan & 1 & 0.191 & 20 & \begin{tikzpicture}[yscale=0.6/5,scale=0.5]
                                  \foreach \i/\y/\c in {3/1/NavyBlue,10/5/Gray} {\draw[line width=10pt,\c] (\i,0)--(\i,\y);};
                                  \draw (0.5,0)--(11.5,0);
                                  \end{tikzpicture} \\[-2.2mm] 
United Kingdom & 1 & 0.144 & 12 & \begin{tikzpicture}[yscale=0.6/5,scale=0.5]
                                  \foreach \i/\y/\c in {6/3/ForestGreen,10/5/Gray,11/1/Maroon} {\draw[line width=10pt,\c] (\i,0)--(\i,\y);};
                                  \draw (0.5,0)--(11.5,0);
                                  \end{tikzpicture} \\[-2.2mm] 
Ireland & 1 & 0.144 & 15 & \begin{tikzpicture}[yscale=0.6/4,scale=0.5]
                                  \foreach \i/\y/\c in {1/3/JungleGreen,4/4/WildStrawberry,8/4/Plum,10/2/Gray,11/3/Maroon} {\draw[line width=10pt,\c] (\i,0)--(\i,\y);};
                                  \draw (0.5,0)--(11.5,0);
                                  \end{tikzpicture} \\[-2.2mm] 
Luxembourg & 1 & 0.116 & 17 & \begin{tikzpicture}[yscale=0.6/4,scale=0.5]
                                  \foreach \i/\y/\c in {2/3/Cerulean,3/4/NavyBlue,7/2/Dandelion,11/3/Maroon} {\draw[line width=10pt,\c] (\i,0)--(\i,\y);};
                                  \draw (0.5,0)--(11.5,0);
                                  \end{tikzpicture} \\[-2.2mm] 
Netherlands & 1 & 0.114 & 7 & \begin{tikzpicture}[yscale=0.6/5,scale=0.5]
                                  \foreach \i/\y/\c in {1/2/JungleGreen,2/1/Cerulean,3/1/NavyBlue,5/3/LimeGreen,8/5/Plum,11/5/Maroon} {\draw[line width=10pt,\c] (\i,0)--(\i,\y);};
                                  \draw (0.5,0)--(11.5,0);
                                  \end{tikzpicture} \\[-2.2mm] 
Poland & 1 & 0.102 & 26 & \begin{tikzpicture}[yscale=0.6/4,scale=0.5]
                                  \foreach \i/\y/\c in {5/2/LimeGreen,10/4/Gray} {\draw[line width=10pt,\c] (\i,0)--(\i,\y);};
                                  \draw (0.5,0)--(11.5,0);
                                  \end{tikzpicture} \\[-2.2mm] 
Germany & 1 & 0.009 & 12 & \begin{tikzpicture}[yscale=0.6/5,scale=0.5]
                                  \foreach \i/\y/\c in {2/2/Cerulean,3/1/NavyBlue,5/3/LimeGreen,6/2/ForestGreen,10/5/Gray,11/2/Maroon} {\draw[line width=10pt,\c] (\i,0)--(\i,\y);};
                                  \draw (0.5,0)--(11.5,0);
                                  \end{tikzpicture} \\[-2.2mm] 
Korea & 2 & 0.275 & 25 & \begin{tikzpicture}[yscale=0.6/5,scale=0.5]
                                  \foreach \i/\y/\c in {5/2/LimeGreen,7/3/Dandelion,10/5/Gray} {\draw[line width=10pt,\c] (\i,0)--(\i,\y);};
                                  \draw (0.5,0)--(11.5,0);
                                  \end{tikzpicture} \\[-2.2mm] 
Austria & 2 & 0.093 & 15 & \begin{tikzpicture}[yscale=0.6/4,scale=0.5]
                                  \foreach \i/\y/\c in {2/2/Cerulean,3/4/NavyBlue,4/4/WildStrawberry,9/1/RedOrange} {\draw[line width=10pt,\c] (\i,0)--(\i,\y);};
                                  \draw (0.5,0)--(11.5,0);
                                  \end{tikzpicture} \\[-2.2mm] 
Estonia & 3 & 0.049 & 28 & \begin{tikzpicture}[yscale=0.6/1,scale=0.5]
                                  \foreach \i/\y/\c in {5/1/LimeGreen} {\draw[line width=10pt,\c] (\i,0)--(\i,\y);};
                                  \draw (0.5,0)--(11.5,0);
                                  \end{tikzpicture} \\[-2.2mm] 
Israel & 5 & 0.116 & 24 & \begin{tikzpicture}[yscale=0.6/1,scale=0.5]
                                  \foreach \i/\y/\c in {8/1/Plum} {\draw[line width=10pt,\c] (\i,0)--(\i,\y);};
                                  \draw (0.5,0)--(11.5,0);
                                  \end{tikzpicture} \\[-2.2mm] 
Czech Republic & 5 & 0.089 & 23 & \begin{tikzpicture}[yscale=0.6/5,scale=0.5]
                                  \foreach \i/\y/\c in {5/5/LimeGreen,10/4/Gray,11/2/Maroon} {\draw[line width=10pt,\c] (\i,0)--(\i,\y);};
                                  \draw (0.5,0)--(11.5,0);
                                  \end{tikzpicture} \\[-2.2mm] 
Slovenia & 6 & 0.091 & 19 & \begin{tikzpicture}[yscale=0.6/3,scale=0.5]
                                  \foreach \i/\y/\c in {4/1/WildStrawberry,5/3/LimeGreen,7/3/Dandelion} {\draw[line width=10pt,\c] (\i,0)--(\i,\y);};
                                  \draw (0.5,0)--(11.5,0);
                                  \end{tikzpicture} \\[-2.2mm] 
France & 6 & 0.066 & 18 & \begin{tikzpicture}[yscale=0.6/2,scale=0.5]
                                  \foreach \i/\y/\c in {6/2/ForestGreen,11/1/Maroon} {\draw[line width=10pt,\c] (\i,0)--(\i,\y);};
                                  \draw (0.5,0)--(11.5,0);
                                  \end{tikzpicture} \\[-2.2mm] 
Mexico & 7 & 0.323 & 35 & \begin{tikzpicture}[yscale=0.6/1,scale=0.5]
                                  \foreach \i/\y/\c in {9/1/RedOrange} {\draw[line width=10pt,\c] (\i,0)--(\i,\y);};
                                  \draw (0.5,0)--(11.5,0);
                                  \end{tikzpicture} \\[-2.2mm] 
Hungary & 8 & 0.058 & 29 & \begin{tikzpicture}[yscale=0.6/1,scale=0.5]
                                  \foreach \i/\y/\c in {10/1/Gray,11/1/Maroon} {\draw[line width=10pt,\c] (\i,0)--(\i,\y);};
                                  \draw (0.5,0)--(11.5,0);
                                  \end{tikzpicture} \\[-2.2mm] 
Russian Federation & 8 & 0.047 & 33 & \begin{tikzpicture}[yscale=0.6/1,scale=0.5]
                                  \foreach \i/\y/\c in {11/1/Maroon} {\draw[line width=10pt,\c] (\i,0)--(\i,\y);};
                                  \draw (0.5,0)--(11.5,0);
                                  \end{tikzpicture} \\[-2.2mm] 
Slovak Republic & 8 & 0.014 & 26 & \begin{tikzpicture}[yscale=0.6/5,scale=0.5]
                                  \foreach \i/\y/\c in {6/4/ForestGreen,10/5/Gray,11/5/Maroon} {\draw[line width=10pt,\c] (\i,0)--(\i,\y);};
                                  \draw (0.5,0)--(11.5,0);
                                  \end{tikzpicture} \\[-2.2mm] 
Turkey & 10 & 0.203 & 36 & \begin{tikzpicture}[yscale=0.6/1,scale=0.5]
                                  \foreach \i/\y/\c in {7/1/Dandelion} {\draw[line width=10pt,\c] (\i,0)--(\i,\y);};
                                  \draw (0.5,0)--(11.5,0);
                                  \end{tikzpicture} \\[-2.2mm] 
Portugal & 10 & 0.12 & 29 & \begin{tikzpicture}[yscale=0.6/5,scale=0.5]
                                  \foreach \i/\y/\c in {1/5/JungleGreen} {\draw[line width=10pt,\c] (\i,0)--(\i,\y);};
                                  \draw (0.5,0)--(11.5,0);
                                  \end{tikzpicture} \\[-2.2mm] 
Italy & 11 & 0.12 & 21 & \begin{tikzpicture}[yscale=0.6/5,scale=0.5]
                                  \foreach \i/\y/\c in {8/1/Plum,11/5/Maroon} {\draw[line width=10pt,\c] (\i,0)--(\i,\y);};
                                  \draw (0.5,0)--(11.5,0);
                                  \end{tikzpicture} \\[-2.2mm] 
Greece & 11 & 0.021 & 34 & \begin{tikzpicture}[yscale=0.6/4,scale=0.5]
                                  \foreach \i/\y/\c in {8/4/Plum,10/3/Gray} {\draw[line width=10pt,\c] (\i,0)--(\i,\y);};
                                  \draw (0.5,0)--(11.5,0);
                                  \end{tikzpicture} \\[-2.2mm] 
Brazil & 12 & 0.209 & 32 & \begin{tikzpicture}[yscale=0.6/5,scale=0.5]
                                  \foreach \i/\y/\c in {4/1/WildStrawberry,9/5/RedOrange} {\draw[line width=10pt,\c] (\i,0)--(\i,\y);};
                                  \draw (0.5,0)--(11.5,0);
                                  \end{tikzpicture} \\[-2.2mm] 
Chile & 21 & 0.169 & 31 & \begin{tikzpicture}[yscale=0.6/2,scale=0.5]
                                  \foreach \i/\y/\c in {9/2/RedOrange,10/1/Gray} {\draw[line width=10pt,\c] (\i,0)--(\i,\y);};
                                  \draw (0.5,0)--(11.5,0);
                                  \end{tikzpicture} \\[-2.2mm] 
\noalign{\smallskip}\hline
\end{tabular}
\end{table}

Optimal weights computed through our optimization tool typically concentrate on some strong dimensions while the other dimensions receive weight zero. On average there are about three non-zero dimensions per country, but only four countries achieve their best position by focusing on only one dimension in 2014 (Mexico, Portugal, Turkey, and the USA, in 2013 this were only Israel and Turkey). This shows that rank-optimal weighting is usually not the same as focusing on a countries best dimension. The largest number of non-zero dimensions in rank-optimal weights is six. Germany and the Netherlands score best with six non-zero dimensions in 2014.

\section{Discussion}

Our method shares some ideas of the ``Benefit of the doubt'' (BOD) method listed in the OECD Handbook \citep{Nardo.Saisana.ea2005HandbookConstructingComposite}, which has recently been applied to the Better Life dataset \citep{Mizobuchi2014MeasuringWorldBetter}. The BOD approach allows country specific weights, thus cross-country comparison based on the same yardstick is not possible. Our approach instead compares common yardsticks (weight vectors) for all countries and picks the one that is best in terms of rank-optimality for a particular country.

Our method of rank optimal weighting extends the toolbox of sensitivity analysis for composite indicators. It enables to find maximal bounds on achievable ranks for each country regarding weighting. The core of our method is the formulation of the optimization problem for rank-optimality in Subsection \ref{subsec:rankmaximization}. This method can be extended or modified in some directions which could also be of interest. 

The method could of course also be used to find the worst possible weights for countries, i.e. weights that bring a country most down in the ranking. 

Some countries achieve their best possible rank by giving a positive weight to only one dimension, which is not in line with the idea of better life as an essentially multi-dimensional concept. To take this into account one could introduce a minimal weight $w^{\min}$ and replace the constraints $w_q\geq 0$ with $w_q\geq w^{\min}$ in the optimization problem. 

Another option would be to define other second order criteria, as already outlined in Footnote \ref{ftnt}. For example, instead of maximizing the distance to the next one in the ranking, one could try to find weight vectors with maximal minimal weight among all those weight vectors which bring the country to its best possible rank. 

The use of composite indicators of multi-dimensional concepts will remain an important tool of data communication, data-driven narratives and policy-making. At the same time awareness increases that choices in the construction can have a critical leverage on the final results. The idea to open up choices in the construction to users as taken by the OECD Better Life Index is valuable to avoid these concerns, but makes a quick assessment of results difficult. Our approach enables a fast and informative first impression on what is maximally feasible through a common yardstick focusing on the performance of a specific country. 

We consider sensitivity analysis and quick exploration of possible rankings of multidimensional datasets the most important purposes for our method, but besides these, one can think of some other applications how countries or individuals can practically use rank-optimal weights for a particular country. As already outlined in the results section, countries could use the rank-optimal weights as target group profile for immigrants, because people of that profile would be comparably satisfied in that country, or even maximally satisfied if the country can make it to the top of the ranking. Individuals who already have a clear picture of which countries to consider as ``best for living'' can probably deduce something about their unconscious preferences by studying the rank-optimal weights of that country. 

Besides multidimensional social indicators the method can also be used for other multidimensional rankings, e.g. for universities or companies. When a university computes its rank-optimal weight in a university ranking it could find out with which profile of importance it outperforms most other universities. This profile could be a good guideline for marketing targeted towards students and faculty.

\paragraph{Acknowledgements.} Jan Lorenz' work financed by the German Research Council (DFG) grant LO2024-1.

\appendix
\section{Appendix}

Tables \ref{tab:2013c} and \ref{tab:2014c} show the rank-optimal weights with continuous weights for the OECD Better Life Index in 2013 and 2014 respectively. The original numbers are available in \citet{Lorenz.Brauer.ea2015RankoptimalWeighting}.

\begin{table}
\centering
\caption{Continuous rank-optimal weights for the OECD Better Life Index 2013.}
\label{tab:2013c}       
\renewcommand{\arraystretch}{0}
\begin{tabular}{m{30mm} >{\centering\bf}m{3mm} >{\centering}m{5mm} >{\centering}m{3mm} m{57mm} }
\hline\noalign{\smallskip}
Country & \rotatebox{90}{Top rank} & \rotatebox{90}{Distance to next} & \rotatebox{90}{Rank equal weights} &\begin{tikzpicture}[scale=0.5] 
  \foreach \i/\y  in {1/Housing,2/Income,3/Jobs,4/Community,5/Education,6/Environment,7/Civic Engagement,8/Health,9/Life Satisfaction,10/Safety,11/Work Life Balance} {\draw (\i,0) node[anchor=south,rotate=90]{\y};}; 
  \end{tikzpicture}  \\ \midrule
United States & 1 & 2.194 & 6 & \begin{tikzpicture}[yscale=0.6/0.872,scale=0.5]
  \foreach \i/\y/\c in {2/0.872/Cerulean,7/0.128/Dandelion} {\draw[line width=10pt,\c] (\i,0)--(\i,\y);};
  \draw (0.5,0)--(11.5,0);
                                  \end{tikzpicture} \\[-2.2mm] 
Switzerland & 1 & 1.297 & 5 & \begin{tikzpicture}[yscale=0.6/0.731,scale=0.5]
  \foreach \i/\y/\c in {2/0.269/Cerulean,9/0.731/RedOrange} {\draw[line width=10pt,\c] (\i,0)--(\i,\y);};
  \draw (0.5,0)--(11.5,0);
                                  \end{tikzpicture} \\[-2.2mm] 
Iceland & 1 & 1.069 & 9 & \begin{tikzpicture}[yscale=0.6/0.763,scale=0.5]
  \foreach \i/\y/\c in {4/0.763/WildStrawberry,9/0.237/RedOrange} {\draw[line width=10pt,\c] (\i,0)--(\i,\y);};
  \draw (0.5,0)--(11.5,0);
                                  \end{tikzpicture} \\[-2.2mm] 
Australia & 1 & 1.063 & 1 & \begin{tikzpicture}[yscale=0.6/0.652,scale=0.5]
  \foreach \i/\y/\c in {1/0.652/JungleGreen,7/0.348/Dandelion} {\draw[line width=10pt,\c] (\i,0)--(\i,\y);};
  \draw (0.5,0)--(11.5,0);
                                  \end{tikzpicture} \\[-2.2mm] 
Finland & 1 & 0.954 & 12 & \begin{tikzpicture}[yscale=0.6/0.882,scale=0.5]
  \foreach \i/\y/\c in {5/0.882/LimeGreen,9/0.118/RedOrange} {\draw[line width=10pt,\c] (\i,0)--(\i,\y);};
  \draw (0.5,0)--(11.5,0);
                                  \end{tikzpicture} \\[-2.2mm] 
Sweden & 1 & 0.762 & 1 & \begin{tikzpicture}[yscale=0.6/0.63,scale=0.5]
  \foreach \i/\y/\c in {5/0.16/LimeGreen,6/0.63/ForestGreen,7/0.128/Dandelion,11/0.082/Maroon} {\draw[line width=10pt,\c] (\i,0)--(\i,\y);};
  \draw (0.5,0)--(11.5,0);
                                  \end{tikzpicture} \\[-2.2mm] 
Norway & 1 & 0.694 & 3 & \begin{tikzpicture}[yscale=0.6/0.318,scale=0.5]
  \foreach \i/\y/\c in {1/0.26/JungleGreen,3/0.318/NavyBlue,6/0.23/ForestGreen,9/0.011/RedOrange,11/0.181/Maroon} {\draw[line width=10pt,\c] (\i,0)--(\i,\y);};
  \draw (0.5,0)--(11.5,0);
                                  \end{tikzpicture} \\[-2.2mm] 
Denmark & 1 & 0.604 & 7 & \begin{tikzpicture}[yscale=0.6/0.824,scale=0.5]
  \foreach \i/\y/\c in {6/0.073/ForestGreen,7/0.102/Dandelion,11/0.824/Maroon} {\draw[line width=10pt,\c] (\i,0)--(\i,\y);};
  \draw (0.5,0)--(11.5,0);
                                  \end{tikzpicture} \\[-2.2mm] 
Japan & 1 & 0.569 & 21 & \begin{tikzpicture}[yscale=0.6/0.521,scale=0.5]
  \foreach \i/\y/\c in {2/0.146/Cerulean,5/0.333/LimeGreen,10/0.521/Gray} {\draw[line width=10pt,\c] (\i,0)--(\i,\y);};
  \draw (0.5,0)--(11.5,0);
                                  \end{tikzpicture} \\[-2.2mm] 
United Kingdom & 1 & 0.519 & 9 & \begin{tikzpicture}[yscale=0.6/0.517,scale=0.5]
  \foreach \i/\y/\c in {2/0.124/Cerulean,4/0.012/WildStrawberry,6/0.517/ForestGreen,10/0.348/Gray} {\draw[line width=10pt,\c] (\i,0)--(\i,\y);};
  \draw (0.5,0)--(11.5,0);
                                  \end{tikzpicture} \\[-2.2mm] 
Ireland & 1 & 0.496 & 14 & \begin{tikzpicture}[yscale=0.6/0.481,scale=0.5]
  \foreach \i/\y/\c in {1/0.323/JungleGreen,4/0.481/WildStrawberry,7/0.058/Dandelion,11/0.138/Maroon} {\draw[line width=10pt,\c] (\i,0)--(\i,\y);};
  \draw (0.5,0)--(11.5,0);
                                  \end{tikzpicture} \\[-2.2mm] 
Luxembourg & 1 & 0.484 & 12 & \begin{tikzpicture}[yscale=0.6/0.383,scale=0.5]
  \foreach \i/\y/\c in {2/0.182/Cerulean,3/0.383/NavyBlue,7/0.199/Dandelion,11/0.236/Maroon} {\draw[line width=10pt,\c] (\i,0)--(\i,\y);};
  \draw (0.5,0)--(11.5,0);
                                  \end{tikzpicture} \\[-2.2mm] 
Canada & 1 & 0.406 & 3 & \begin{tikzpicture}[yscale=0.6/0.72,scale=0.5]
  \foreach \i/\y/\c in {1/0.006/JungleGreen,2/0.086/Cerulean,9/0.123/RedOrange,10/0.72/Gray,11/0.065/Maroon} {\draw[line width=10pt,\c] (\i,0)--(\i,\y);};
  \draw (0.5,0)--(11.5,0);
                                  \end{tikzpicture} \\[-2.2mm] 
Netherlands & 1 & 0.364 & 7 & \begin{tikzpicture}[yscale=0.6/0.433,scale=0.5]
  \foreach \i/\y/\c in {1/0.136/JungleGreen,2/0.16/Cerulean,3/0.143/NavyBlue,4/0.128/WildStrawberry,11/0.433/Maroon} {\draw[line width=10pt,\c] (\i,0)--(\i,\y);};
  \draw (0.5,0)--(11.5,0);
                                  \end{tikzpicture} \\[-2.2mm] 
Belgium & 1 & 0.102 & 16 & \begin{tikzpicture}[yscale=0.6/0.485,scale=0.5]
  \foreach \i/\y/\c in {2/0.284/Cerulean,5/0.105/LimeGreen,7/0.126/Dandelion,11/0.485/Maroon} {\draw[line width=10pt,\c] (\i,0)--(\i,\y);};
  \draw (0.5,0)--(11.5,0);
                                  \end{tikzpicture} \\[-2.2mm] 
New Zealand & 1 & 0.088 & 9 & \begin{tikzpicture}[yscale=0.6/0.616,scale=0.5]
  \foreach \i/\y/\c in {6/0.143/ForestGreen,8/0.616/Plum,10/0.108/Gray,11/0.132/Maroon} {\draw[line width=10pt,\c] (\i,0)--(\i,\y);};
  \draw (0.5,0)--(11.5,0);
                                  \end{tikzpicture} \\[-2.2mm] 
Spain & 1 & 0.082 & 20 & \begin{tikzpicture}[yscale=0.6/0.459,scale=0.5]
  \foreach \i/\y/\c in {8/0.459/Plum,10/0.187/Gray,11/0.354/Maroon} {\draw[line width=10pt,\c] (\i,0)--(\i,\y);};
  \draw (0.5,0)--(11.5,0);
                                  \end{tikzpicture} \\[-2.2mm] 
Poland & 2 & 0.101 & 25 & \begin{tikzpicture}[yscale=0.6/0.785,scale=0.5]
  \foreach \i/\y/\c in {5/0.068/LimeGreen,10/0.785/Gray,11/0.147/Maroon} {\draw[line width=10pt,\c] (\i,0)--(\i,\y);};
  \draw (0.5,0)--(11.5,0);
                                  \end{tikzpicture} \\[-2.2mm] 
Korea & 2 & 0.02 & 26 & \begin{tikzpicture}[yscale=0.6/0.487,scale=0.5]
  \foreach \i/\y/\c in {5/0.24/LimeGreen,7/0.273/Dandelion,10/0.487/Gray} {\draw[line width=10pt,\c] (\i,0)--(\i,\y);};
  \draw (0.5,0)--(11.5,0);
                                  \end{tikzpicture} \\[-2.2mm] 
Austria & 2 & 0.008 & 14 & \begin{tikzpicture}[yscale=0.6/0.221,scale=0.5]
  \foreach \i/\y/\c in {2/0.11/Cerulean,3/0.221/NavyBlue,4/0.187/WildStrawberry,7/0.197/Dandelion,9/0.125/RedOrange,10/0.16/Gray} {\draw[line width=10pt,\c] (\i,0)--(\i,\y);};
  \draw (0.5,0)--(11.5,0);
                                  \end{tikzpicture} \\[-2.2mm] 
Germany & 3 & 0.118 & 17 & \begin{tikzpicture}[yscale=0.6/0.477,scale=0.5]
  \foreach \i/\y/\c in {2/0.092/Cerulean,5/0.477/LimeGreen,6/0.412/ForestGreen,11/0.019/Maroon} {\draw[line width=10pt,\c] (\i,0)--(\i,\y);};
  \draw (0.5,0)--(11.5,0);
                                  \end{tikzpicture} \\[-2.2mm] 
Israel & 5 & 0.128 & 24 & \begin{tikzpicture}[yscale=0.6/1,scale=0.5]
  \foreach \i/\y/\c in {8/1/Plum} {\draw[line width=10pt,\c] (\i,0)--(\i,\y);};
  \draw (0.5,0)--(11.5,0);
                                  \end{tikzpicture} \\[-2.2mm] 
Slovenia & 5 & 0.011 & 19 & \begin{tikzpicture}[yscale=0.6/0.622,scale=0.5]
  \foreach \i/\y/\c in {5/0.622/LimeGreen,7/0.163/Dandelion,10/0.019/Gray,11/0.195/Maroon} {\draw[line width=10pt,\c] (\i,0)--(\i,\y);};
  \draw (0.5,0)--(11.5,0);
                                  \end{tikzpicture} \\[-2.2mm] 
Russian Federation & 7 & 0.019 & 32 & \begin{tikzpicture}[yscale=0.6/0.998,scale=0.5]
  \foreach \i/\y/\c in {3/0.002/NavyBlue,11/0.998/Maroon} {\draw[line width=10pt,\c] (\i,0)--(\i,\y);};
  \draw (0.5,0)--(11.5,0);
                                  \end{tikzpicture} \\[-2.2mm] 
France & 8 & 0.094 & 18 & \begin{tikzpicture}[yscale=0.6/0.403,scale=0.5]
  \foreach \i/\y/\c in {1/0.152/JungleGreen,2/0.403/Cerulean,4/0.203/WildStrawberry,8/0.082/Plum,11/0.16/Maroon} {\draw[line width=10pt,\c] (\i,0)--(\i,\y);};
  \draw (0.5,0)--(11.5,0);
                                  \end{tikzpicture} \\[-2.2mm] 
Portugal & 8 & 0.066 & 28 & \begin{tikzpicture}[yscale=0.6/0.893,scale=0.5]
  \foreach \i/\y/\c in {1/0.893/JungleGreen,6/0.107/ForestGreen} {\draw[line width=10pt,\c] (\i,0)--(\i,\y);};
  \draw (0.5,0)--(11.5,0);
                                  \end{tikzpicture} \\[-2.2mm] 
Hungary & 8 & 0.027 & 29 & \begin{tikzpicture}[yscale=0.6/0.536,scale=0.5]
  \foreach \i/\y/\c in {10/0.464/Gray,11/0.536/Maroon} {\draw[line width=10pt,\c] (\i,0)--(\i,\y);};
  \draw (0.5,0)--(11.5,0);
                                  \end{tikzpicture} \\[-2.2mm] 
Estonia & 8 & 0.008 & 31 & \begin{tikzpicture}[yscale=0.6/0.803,scale=0.5]
  \foreach \i/\y/\c in {5/0.803/LimeGreen,6/0.072/ForestGreen,11/0.125/Maroon} {\draw[line width=10pt,\c] (\i,0)--(\i,\y);};
  \draw (0.5,0)--(11.5,0);
                                  \end{tikzpicture} \\[-2.2mm] 
Italy & 9 & 0.026 & 22 & \begin{tikzpicture}[yscale=0.6/0.622,scale=0.5]
  \foreach \i/\y/\c in {2/0.147/Cerulean,7/0.194/Dandelion,8/0.038/Plum,11/0.622/Maroon} {\draw[line width=10pt,\c] (\i,0)--(\i,\y);};
  \draw (0.5,0)--(11.5,0);
                                  \end{tikzpicture} \\[-2.2mm] 
Slovak Republic & 9 & 0.01 & 26 & \begin{tikzpicture}[yscale=0.6/0.596,scale=0.5]
  \foreach \i/\y/\c in {6/0.102/ForestGreen,10/0.596/Gray,11/0.302/Maroon} {\draw[line width=10pt,\c] (\i,0)--(\i,\y);};
  \draw (0.5,0)--(11.5,0);
                                  \end{tikzpicture} \\[-2.2mm] 
Mexico & 10 & 0.323 & 35 & \begin{tikzpicture}[yscale=0.6/1,scale=0.5]
  \foreach \i/\y/\c in {9/1/RedOrange} {\draw[line width=10pt,\c] (\i,0)--(\i,\y);};
  \draw (0.5,0)--(11.5,0);
                                  \end{tikzpicture} \\[-2.2mm] 
Czech Republic & 10 & 0.047 & 22 & \begin{tikzpicture}[yscale=0.6/0.615,scale=0.5]
  \foreach \i/\y/\c in {5/0.385/LimeGreen,10/0.615/Gray} {\draw[line width=10pt,\c] (\i,0)--(\i,\y);};
  \draw (0.5,0)--(11.5,0);
                                  \end{tikzpicture} \\[-2.2mm] 
Turkey & 11 & 0.114 & 36 & \begin{tikzpicture}[yscale=0.6/1,scale=0.5]
  \foreach \i/\y/\c in {7/1/Dandelion} {\draw[line width=10pt,\c] (\i,0)--(\i,\y);};
  \draw (0.5,0)--(11.5,0);
                                  \end{tikzpicture} \\[-2.2mm] 
Greece & 12 & 0.04 & 29 & \begin{tikzpicture}[yscale=0.6/0.571,scale=0.5]
  \foreach \i/\y/\c in {8/0.571/Plum,10/0.424/Gray,11/0.004/Maroon} {\draw[line width=10pt,\c] (\i,0)--(\i,\y);};
  \draw (0.5,0)--(11.5,0);
                                  \end{tikzpicture} \\[-2.2mm] 
Brazil & 18 & 0.153 & 33 & \begin{tikzpicture}[yscale=0.6/0.583,scale=0.5]
  \foreach \i/\y/\c in {7/0.387/Dandelion,9/0.583/RedOrange,11/0.03/Maroon} {\draw[line width=10pt,\c] (\i,0)--(\i,\y);};
  \draw (0.5,0)--(11.5,0);
                                  \end{tikzpicture} \\[-2.2mm] 
Chile & 21 & 0.081 & 34 & \begin{tikzpicture}[yscale=0.6/0.815,scale=0.5]
  \foreach \i/\y/\c in {8/0.017/Plum,9/0.815/RedOrange,10/0.167/Gray} {\draw[line width=10pt,\c] (\i,0)--(\i,\y);};
  \draw (0.5,0)--(11.5,0);
                                  \end{tikzpicture} \\[-2.2mm] 
\noalign{\smallskip}\hline
\end{tabular}
\end{table}

\begin{table}
\centering
\caption{Continuous rank-optimal weights for the OECD Better Life Index 2014.}
\label{tab:2014c}       
\renewcommand{\arraystretch}{0}
\begin{tabular}{m{30mm} >{\centering\bf}m{3mm} >{\centering}m{5mm} >{\centering}m{3mm} m{57mm} }
\hline\noalign{\smallskip}
Country & \rotatebox{90}{Top rank} & \rotatebox{90}{Distance to next} & \rotatebox{90}{Rank equal weights} &\begin{tikzpicture}[scale=0.5] 
  \foreach \i/\y  in {1/Housing,2/Income,3/Jobs,4/Community,5/Education,6/Environment,7/Civic Engagement,8/Health,9/Life Satisfaction,10/Safety,11/Work Life Balance} {\draw (\i,0) node[anchor=south,rotate=90]{\y};}; 
  \end{tikzpicture}  \\ \midrule
United States & 1 & 2.739 & 7 & \begin{tikzpicture}[yscale=0.6/1,scale=0.5]
  \foreach \i/\y/\c in {2/1/Cerulean} {\draw[line width=10pt,\c] (\i,0)--(\i,\y);};
  \draw (0.5,0)--(11.5,0);
                                  \end{tikzpicture} \\[-2.2mm] 
Finland & 1 & 1.15 & 7 & \begin{tikzpicture}[yscale=0.6/0.916,scale=0.5]
  \foreach \i/\y/\c in {5/0.916/LimeGreen,9/0.084/RedOrange} {\draw[line width=10pt,\c] (\i,0)--(\i,\y);};
  \draw (0.5,0)--(11.5,0);
                                  \end{tikzpicture} \\[-2.2mm] 
Australia & 1 & 1.093 & 1 & \begin{tikzpicture}[yscale=0.6/0.387,scale=0.5]
  \foreach \i/\y/\c in {1/0.28/JungleGreen,7/0.334/Dandelion,10/0.387/Gray} {\draw[line width=10pt,\c] (\i,0)--(\i,\y);};
  \draw (0.5,0)--(11.5,0);
                                  \end{tikzpicture} \\[-2.2mm] 
Switzerland & 1 & 1.066 & 4 & \begin{tikzpicture}[yscale=0.6/0.533,scale=0.5]
  \foreach \i/\y/\c in {2/0.244/Cerulean,3/0.224/NavyBlue,9/0.533/RedOrange} {\draw[line width=10pt,\c] (\i,0)--(\i,\y);};
  \draw (0.5,0)--(11.5,0);
                                  \end{tikzpicture} \\[-2.2mm] 
Denmark & 1 & 1.05 & 4 & \begin{tikzpicture}[yscale=0.6/0.555,scale=0.5]
  \foreach \i/\y/\c in {4/0.555/WildStrawberry,9/0.014/RedOrange,11/0.431/Maroon} {\draw[line width=10pt,\c] (\i,0)--(\i,\y);};
  \draw (0.5,0)--(11.5,0);
                                  \end{tikzpicture} \\[-2.2mm] 
Sweden & 1 & 0.819 & 2 & \begin{tikzpicture}[yscale=0.6/0.982,scale=0.5]
  \foreach \i/\y/\c in {6/0.982/ForestGreen,11/0.018/Maroon} {\draw[line width=10pt,\c] (\i,0)--(\i,\y);};
  \draw (0.5,0)--(11.5,0);
                                  \end{tikzpicture} \\[-2.2mm] 
Norway & 1 & 0.771 & 2 & \begin{tikzpicture}[yscale=0.6/0.424,scale=0.5]
  \foreach \i/\y/\c in {1/0.424/JungleGreen,3/0.279/NavyBlue,6/0.148/ForestGreen,11/0.149/Maroon} {\draw[line width=10pt,\c] (\i,0)--(\i,\y);};
  \draw (0.5,0)--(11.5,0);
                                  \end{tikzpicture} \\[-2.2mm] 
Spain & 1 & 0.429 & 21 & \begin{tikzpicture}[yscale=0.6/0.521,scale=0.5]
  \foreach \i/\y/\c in {8/0.521/Plum,11/0.479/Maroon} {\draw[line width=10pt,\c] (\i,0)--(\i,\y);};
  \draw (0.5,0)--(11.5,0);
                                  \end{tikzpicture} \\[-2.2mm] 
New Zealand & 1 & 0.396 & 10 & \begin{tikzpicture}[yscale=0.6/0.488,scale=0.5]
  \foreach \i/\y/\c in {4/0.488/WildStrawberry,7/0.097/Dandelion,8/0.336/Plum,11/0.079/Maroon} {\draw[line width=10pt,\c] (\i,0)--(\i,\y);};
  \draw (0.5,0)--(11.5,0);
                                  \end{tikzpicture} \\[-2.2mm] 
Canada & 1 & 0.383 & 4 & \begin{tikzpicture}[yscale=0.6/0.313,scale=0.5]
  \foreach \i/\y/\c in {1/0.283/JungleGreen,2/0.093/Cerulean,9/0.311/RedOrange,10/0.313/Gray} {\draw[line width=10pt,\c] (\i,0)--(\i,\y);};
  \draw (0.5,0)--(11.5,0);
                                  \end{tikzpicture} \\[-2.2mm] 
Belgium & 1 & 0.345 & 12 & \begin{tikzpicture}[yscale=0.6/0.403,scale=0.5]
  \foreach \i/\y/\c in {1/0.317/JungleGreen,2/0.223/Cerulean,7/0.057/Dandelion,11/0.403/Maroon} {\draw[line width=10pt,\c] (\i,0)--(\i,\y);};
  \draw (0.5,0)--(11.5,0);
                                  \end{tikzpicture} \\[-2.2mm] 
Iceland & 1 & 0.256 & 10 & \begin{tikzpicture}[yscale=0.6/0.516,scale=0.5]
  \foreach \i/\y/\c in {3/0.24/NavyBlue,4/0.516/WildStrawberry,8/0.018/Plum,10/0.226/Gray} {\draw[line width=10pt,\c] (\i,0)--(\i,\y);};
  \draw (0.5,0)--(11.5,0);
                                  \end{tikzpicture} \\[-2.2mm] 
Japan & 1 & 0.229 & 20 & \begin{tikzpicture}[yscale=0.6/0.979,scale=0.5]
  \foreach \i/\y/\c in {2/0.021/Cerulean,10/0.979/Gray} {\draw[line width=10pt,\c] (\i,0)--(\i,\y);};
  \draw (0.5,0)--(11.5,0);
                                  \end{tikzpicture} \\[-2.2mm] 
United Kingdom & 1 & 0.177 & 12 & \begin{tikzpicture}[yscale=0.6/0.553,scale=0.5]
  \foreach \i/\y/\c in {4/0.014/WildStrawberry,6/0.335/ForestGreen,10/0.553/Gray,11/0.098/Maroon} {\draw[line width=10pt,\c] (\i,0)--(\i,\y);};
  \draw (0.5,0)--(11.5,0);
                                  \end{tikzpicture} \\[-2.2mm] 
Ireland & 1 & 0.153 & 15 & \begin{tikzpicture}[yscale=0.6/0.261,scale=0.5]
  \foreach \i/\y/\c in {1/0.199/JungleGreen,4/0.26/WildStrawberry,8/0.261/Plum,10/0.095/Gray,11/0.186/Maroon} {\draw[line width=10pt,\c] (\i,0)--(\i,\y);};
  \draw (0.5,0)--(11.5,0);
                                  \end{tikzpicture} \\[-2.2mm] 
Luxembourg & 1 & 0.149 & 17 & \begin{tikzpicture}[yscale=0.6/0.332,scale=0.5]
  \foreach \i/\y/\c in {2/0.245/Cerulean,3/0.332/NavyBlue,7/0.166/Dandelion,11/0.257/Maroon} {\draw[line width=10pt,\c] (\i,0)--(\i,\y);};
  \draw (0.5,0)--(11.5,0);
                                  \end{tikzpicture} \\[-2.2mm] 
Poland & 1 & 0.148 & 26 & \begin{tikzpicture}[yscale=0.6/0.665,scale=0.5]
  \foreach \i/\y/\c in {5/0.28/LimeGreen,7/0.055/Dandelion,10/0.665/Gray} {\draw[line width=10pt,\c] (\i,0)--(\i,\y);};
  \draw (0.5,0)--(11.5,0);
                                  \end{tikzpicture} \\[-2.2mm] 
Netherlands & 1 & 0.138 & 7 & \begin{tikzpicture}[yscale=0.6/0.396,scale=0.5]
  \foreach \i/\y/\c in {1/0.075/JungleGreen,2/0.069/Cerulean,3/0.039/NavyBlue,5/0.036/LimeGreen,8/0.396/Plum,11/0.386/Maroon} {\draw[line width=10pt,\c] (\i,0)--(\i,\y);};
  \draw (0.5,0)--(11.5,0);
                                  \end{tikzpicture} \\[-2.2mm] 
Germany & 1 & 0.027 & 12 & \begin{tikzpicture}[yscale=0.6/0.307,scale=0.5]
  \foreach \i/\y/\c in {2/0.141/Cerulean,3/0.042/NavyBlue,5/0.157/LimeGreen,6/0.241/ForestGreen,10/0.307/Gray,11/0.112/Maroon} {\draw[line width=10pt,\c] (\i,0)--(\i,\y);};
  \draw (0.5,0)--(11.5,0);
                                  \end{tikzpicture} \\[-2.2mm] 
Korea & 2 & 0.341 & 25 & \begin{tikzpicture}[yscale=0.6/0.528,scale=0.5]
  \foreach \i/\y/\c in {5/0.194/LimeGreen,7/0.277/Dandelion,10/0.528/Gray} {\draw[line width=10pt,\c] (\i,0)--(\i,\y);};
  \draw (0.5,0)--(11.5,0);
                                  \end{tikzpicture} \\[-2.2mm] 
Austria & 2 & 0.118 & 15 & \begin{tikzpicture}[yscale=0.6/0.336,scale=0.5]
  \foreach \i/\y/\c in {2/0.18/Cerulean,3/0.315/NavyBlue,4/0.336/WildStrawberry,7/0.039/Dandelion,9/0.13/RedOrange} {\draw[line width=10pt,\c] (\i,0)--(\i,\y);};
  \draw (0.5,0)--(11.5,0);
                                  \end{tikzpicture} \\[-2.2mm] 
Estonia & 3 & 0.064 & 28 & \begin{tikzpicture}[yscale=0.6/0.995,scale=0.5]
  \foreach \i/\y/\c in {5/0.995/LimeGreen,6/0.005/ForestGreen} {\draw[line width=10pt,\c] (\i,0)--(\i,\y);};
  \draw (0.5,0)--(11.5,0);
                                  \end{tikzpicture} \\[-2.2mm] 
Czech Republic & 4 & 0.033 & 23 & \begin{tikzpicture}[yscale=0.6/0.374,scale=0.5]
  \foreach \i/\y/\c in {5/0.374/LimeGreen,6/0.058/ForestGreen,10/0.37/Gray,11/0.198/Maroon} {\draw[line width=10pt,\c] (\i,0)--(\i,\y);};
  \draw (0.5,0)--(11.5,0);
                                  \end{tikzpicture} \\[-2.2mm] 
Israel & 5 & 0.116 & 24 & \begin{tikzpicture}[yscale=0.6/0.989,scale=0.5]
  \foreach \i/\y/\c in {2/0.011/Cerulean,8/0.989/Plum} {\draw[line width=10pt,\c] (\i,0)--(\i,\y);};
  \draw (0.5,0)--(11.5,0);
                                  \end{tikzpicture} \\[-2.2mm] 
France & 6 & 0.115 & 18 & \begin{tikzpicture}[yscale=0.6/0.639,scale=0.5]
  \foreach \i/\y/\c in {2/0.032/Cerulean,6/0.639/ForestGreen,11/0.329/Maroon} {\draw[line width=10pt,\c] (\i,0)--(\i,\y);};
  \draw (0.5,0)--(11.5,0);
                                  \end{tikzpicture} \\[-2.2mm] 
Slovenia & 6 & 0.115 & 19 & \begin{tikzpicture}[yscale=0.6/0.471,scale=0.5]
  \foreach \i/\y/\c in {4/0.134/WildStrawberry,5/0.471/LimeGreen,7/0.395/Dandelion} {\draw[line width=10pt,\c] (\i,0)--(\i,\y);};
  \draw (0.5,0)--(11.5,0);
                                  \end{tikzpicture} \\[-2.2mm] 
Mexico & 7 & 0.323 & 35 & \begin{tikzpicture}[yscale=0.6/1,scale=0.5]
  \foreach \i/\y/\c in {9/1/RedOrange} {\draw[line width=10pt,\c] (\i,0)--(\i,\y);};
  \draw (0.5,0)--(11.5,0);
                                  \end{tikzpicture} \\[-2.2mm] 
Hungary & 8 & 0.067 & 29 & \begin{tikzpicture}[yscale=0.6/0.497,scale=0.5]
  \foreach \i/\y/\c in {7/0.01/Dandelion,10/0.497/Gray,11/0.493/Maroon} {\draw[line width=10pt,\c] (\i,0)--(\i,\y);};
  \draw (0.5,0)--(11.5,0);
                                  \end{tikzpicture} \\[-2.2mm] 
Russian Federation & 8 & 0.066 & 33 & \begin{tikzpicture}[yscale=0.6/0.785,scale=0.5]
  \foreach \i/\y/\c in {3/0.111/NavyBlue,5/0.104/LimeGreen,11/0.785/Maroon} {\draw[line width=10pt,\c] (\i,0)--(\i,\y);};
  \draw (0.5,0)--(11.5,0);
                                  \end{tikzpicture} \\[-2.2mm] 
Slovak Republic & 8 & 0.015 & 26 & \begin{tikzpicture}[yscale=0.6/0.361,scale=0.5]
  \foreach \i/\y/\c in {6/0.29/ForestGreen,10/0.349/Gray,11/0.361/Maroon} {\draw[line width=10pt,\c] (\i,0)--(\i,\y);};
  \draw (0.5,0)--(11.5,0);
                                  \end{tikzpicture} \\[-2.2mm] 
Turkey & 10 & 0.203 & 36 & \begin{tikzpicture}[yscale=0.6/1,scale=0.5]
  \foreach \i/\y/\c in {7/1/Dandelion} {\draw[line width=10pt,\c] (\i,0)--(\i,\y);};
  \draw (0.5,0)--(11.5,0);
                                  \end{tikzpicture} \\[-2.2mm] 
Portugal & 10 & 0.12 & 29 & \begin{tikzpicture}[yscale=0.6/1,scale=0.5]
  \foreach \i/\y/\c in {1/1/JungleGreen} {\draw[line width=10pt,\c] (\i,0)--(\i,\y);};
  \draw (0.5,0)--(11.5,0);
                                  \end{tikzpicture} \\[-2.2mm] 
Italy & 10 & 0.026 & 21 & \begin{tikzpicture}[yscale=0.6/0.596,scale=0.5]
  \foreach \i/\y/\c in {7/0.13/Dandelion,8/0.274/Plum,11/0.596/Maroon} {\draw[line width=10pt,\c] (\i,0)--(\i,\y);};
  \draw (0.5,0)--(11.5,0);
                                  \end{tikzpicture} \\[-2.2mm] 
Greece & 11 & 0.021 & 34 & \begin{tikzpicture}[yscale=0.6/0.597,scale=0.5]
  \foreach \i/\y/\c in {8/0.597/Plum,10/0.403/Gray} {\draw[line width=10pt,\c] (\i,0)--(\i,\y);};
  \draw (0.5,0)--(11.5,0);
                                  \end{tikzpicture} \\[-2.2mm] 
Brazil & 12 & 0.229 & 32 & \begin{tikzpicture}[yscale=0.6/0.863,scale=0.5]
  \foreach \i/\y/\c in {4/0.137/WildStrawberry,9/0.863/RedOrange} {\draw[line width=10pt,\c] (\i,0)--(\i,\y);};
  \draw (0.5,0)--(11.5,0);
                                  \end{tikzpicture} \\[-2.2mm] 
Chile & 21 & 0.181 & 31 & \begin{tikzpicture}[yscale=0.6/0.669,scale=0.5]
  \foreach \i/\y/\c in {9/0.669/RedOrange,10/0.331/Gray} {\draw[line width=10pt,\c] (\i,0)--(\i,\y);};
  \draw (0.5,0)--(11.5,0);
                                  \end{tikzpicture} \\[-2.2mm] 
\noalign{\smallskip}\hline
\end{tabular}
\end{table}

\bibliographystyle{spbasic}

\end{document}